\documentclass[a4paper,10pt]{article}
\usepackage{geometry}
\usepackage[numbers,square,sort&compress]{natbib}
\usepackage{graphicx,subcaption}
\usepackage{amssymb,amsmath,amsthm}
\usepackage{bm}
\usepackage{lineno}
\usepackage{setspace}
\usepackage{booktabs,multirow}
\usepackage{authblk}
\usepackage{float}
\usepackage[flushleft]{threeparttable}
\usepackage{enumitem}
\usepackage{caption}
\allowdisplaybreaks
\usepackage{color}
\usepackage[colorlinks,citecolor=blue,urlcolor=blue]{hyperref}

\newtheorem{theorem}{Theorem}[section]
\newtheorem{lemma}{Lemma}[section]

\newtheorem{corollary}{Corollary}[section]
\newtheorem{remark}{Remark}[section]

\numberwithin{equation}{section}

\title{Shock solutions of a 3D spherical piston problem for the compressible Euler equations}
\author[]{Pengfei Zhou}
\author[]{Huimin Yu}
\author[]{Xiaomin Zhang \thanks{Corresponding author: zxm15924687@163.com}}
\affil[]{School of Mathematics and Statistics, Shandong Normal University, Jinan 250014, China.}
\date{}
\begin{document}
	
	\maketitle

		\begin{flushleft}
			\textbf{Abstract:} This paper establishes the global well-posedness of shock solutions for a 3D spherical piston problem to the compressible Euler equations. The main mathematical challenges are the nonlinear coupling between the free shock boundary and the moving piston, alongside potential gradient catastrophes induced by spherical geometric source terms.
			To overcome these difficulties, we reformulate the problem using Riemann invariants and introduce a tailored directional derivative along the unperturbed piston trajectory. Then we can extend the local solution to construct a global piecewise smooth solution including a single shock under some moderate assumptions.
		\end{flushleft} 
		\begin{flushleft}
		\end{flushleft}
		\begin{flushleft}
		\end{flushleft}
	\section{Introduction}
	The compressible Euler equations constitute a canonical system of hyperbolic conservation laws for modeling inviscid compressible flow. This model captures the fundamental nonlinear wave phenomena\textemdash including shock waves, rarefaction waves and supersonic-subsonic transitions in aerodynamics, propulsion systems and astrophysics, see \cite{1,2,3,4,5,6} for example.
	In this paper, we consider the following 3D isentropic compressible Euler equations:
	\begin{equation}\label{ce}
		\left\{
		\begin{aligned}
			&\partial_t\rho+\nabla\cdot(\rho\mathbf{u})=0,\\
			&\partial_t(\rho\mathbf{u})+\nabla\cdot(\rho\mathbf{u}\otimes\mathbf{u})+\nabla p(\rho)=0,\\
		\end{aligned}
		\right.
	\end{equation}
	with $t\in\mathbb{R}^+$, $\mathbf{x}=(x_1,x_2,x_3)\in\mathbb{R}^3$, where $\rho$ is the density of fluid, $\mathbf{u}=(u_1,u_2,u_3)^\text{T}$ is the velocity vector field. For polytropic gas, the pressure satisfies $p(\rho)=\rho^\gamma$ with the adiabatic exponent $\gamma\geqslant1$.
	
	To analyze the spherical flow, we introduce the spherical coordinates transformation:
	\begin{equation*}
		\rho(\mathbf{x},t)=\rho(r,t),\quad\mathbf{u}(\mathbf{x},t)=\frac{\mathbf{x}}{r}\cdot u(r,t),
	\end{equation*}
	with the radius $r=\sqrt{x_1^2+x_2^2+x_3^2}$, then system (\ref{ce}) becomes
	\begin{equation}\label{e1}
		\left\{
		\begin{aligned}
			&\partial_t\rho+\partial_r(\rho u)+\frac{2}{r}\rho u=0,\\
			&\partial_t(\rho u)+\partial_r\left(\rho u^2+p\right)+\frac{2}{r}\rho u^2=0,\\
		\end{aligned}
		\right.
	\end{equation}
	 which serves as a simplified model for investigating the existence, regularity, and structural properties of solutions to multi-dimensional inviscid compressible fluid flows. The local existence and uniqueness of classical solutions for \eqref{e1} with smooth initial data away from vacuum was established by Majda \cite{7}.  For small perturbations of a constant equilibrium state, Klainerman \cite{8} developed energy methods and symmetrization techniques to prove the global existence to smooth solutions. The core mechanism by which small initial data suppresses nonlinear blowup and maintains long-time regularity was further elucidated by Christodoulou \cite{9}.
	
	The piston problem is a basic model in gas dynamics \cite{1,10}. Under the small perturbation framework, some significant progress had been made. For the one-dimensional case, Li and Zhao \cite{11} obtained the stability of the shock wave by the method of characteristics under sufficiently small weighted $C^1$ perturbations. Ding \cite{12} established the global
	existence of entropy solutions under the assumptions that the total variations
	of the initial data and the piston velocity are both sufficiently small via a modified Glimm scheme. Zhang and Yu \cite{13} showed the global well-posedness of the shock solution to the non-isentropic compressible Euler equations in $C^1$ space. Hu, Li and Zhang \cite{14} studied the one-dimensional piston problem with large velocity variations, and proved the global existence of piecewise smooth shock solutions under the vanishing initial density assumption. This work broke through the classical small perturbation framework for piston problems and reveals an intrinsic stabilization mechanism in the high Mach number regime.
	
	For multi-dimensional case, \cite{15} derived unsteady self-similar shock solutions for the uniformly expanding spherical piston, constituting a cornerstone result for investigating the spherical piston problem. \cite{16} concerned a special initial boundary value problem of multi-dimensional unsteady potential flow equation. It showed the stability of the structure of shock front in symmetric case under small perturbation, which became the foundation for future research. Subsequently, \cite{17,18,19,20} developed the existence theory of BV solutions for several piston problems of compressible flow. And the global well-posedness of piecewise-smooth solutions remains open for Euler equations. The relevant result for the relativistic Euler equation can be found in \cite{21}.
	
     It is worth noting that generalizing the large velocity variation condition proposed in \cite{14} to the three-dimensional spherical piston problem will meet a series of difficulties that are completely absent in the one-dimensional model. The spherical geometry effect inherently amplifies flow gradients and raises the risk of gradient catastrophe, rendering the study of the three-dimensional spherically symmetric problem infeasible within the large velocity variation framework. Accordingly, under the assumption of vanishing initial density, we introduce radial decay weights $r^{-\delta_0}$ on the perturbations in Theorem \ref{th1}, thereby ensuring the global well-posedness of piecewise-smooth solutions.
     
     Consider a physical scenario in which a spherical piston expands into a quiescent gas, thereby generating a forward shock wave. That is we assume that the gas initially in a uniform static state
	\begin{equation}\label{ini}
		\left(\rho,u\right)\left(r,0\right)=\left(\rho_0,0\right).
	\end{equation}
	Then the fluid is driven by a moving piston $t=b(r)$, a natural boundary condition on piston is that the fluid velocity is equal to the piston
	velocity \cite{1}, i.e.
	\begin{equation}\label{e13}
		t=b(r):\ u\left(r,b(r)\right)=\frac{1}{b'(r)}.
	\end{equation}
	 And the expanding motion generates a forward shock wave $t=s(r)$. Across this shock, the piecewise smooth solution satisfies the Rankine–Hugoniot conditions
	\begin{equation}\label{e14}
		t=s(r):\\
		\left\{
		\begin{aligned}
			&\rho u=\frac{1}{s'(r)}{(\rho-\rho_0)},\\
			&\rho u(u-\frac{1}{s'(r)})+\rho^{\gamma}=\rho_0^{\gamma}.\\
		\end{aligned}
		\right.
	\end{equation}
	
     Due to the time-dependence of the expanding piston boundary and the geometrical effect induced by the spherically‑symmetric geometry, steady-state solutions fail to capture the global dynamic expansion process, which physically motivates us to use self-similar solutions as the background states. As established in \cite{16}, if the piston expands at a constant velocity $b_0>0$ and the initial flow is uniformly quiescent, i.e. \eqref{ini}, the piston problem admits a global in time shock front solution with constant shock speed $\nu$, which is self-similar. We introduce the  variable $\eta=t/r$ and denote $s_0=1/\nu$, then system \eqref{e1} together with boundary conditions \eqref{e13} and \eqref{e14} yields
    \begin{equation}\label{e21}
    	\left\{
    	\begin{aligned}
    		&\left(1-\eta \hat{u}\right)\hat{\rho}_{\eta}-\eta\hat{\rho} \hat{u}_{\eta}+2\hat{\rho} \hat{u}=0,\\
    		&\left(1-\eta \hat{u}\right)\hat{u}_{\eta}-\gamma\eta\hat{\rho}^{\gamma-2}\hat{\rho}_{\eta}=0,\\
    		&\hat{u}=b_0\quad \text{on}\ \eta=\frac{1}{b_0},\\
    		&\hat{\rho} \hat{u}=\nu(\hat{\rho}-\rho_0)\quad \text{on}\ \eta=s_0,\\
    		&\hat{\rho} \hat{u}(\hat{u}-\nu)+\hat{\rho}^{\gamma}=\rho_0^{\gamma}\quad \text{on}\ \eta=s_0,\\
    	\end{aligned}
    	\right.
    \end{equation}
	The existence and uniqueness of solution $(\hat{\rho},\hat{u})(\eta)$ on $\left[s_0,1/b_0\right]$ to the piston problem \eqref{e21} can be obtained by using the method similar to that in \cite{16}.
	
	 In this paper, we prove the global well-posedness of shock front solutions to the piston problem in a small neighbourhood of this self-similar solution $(\hat{\rho},\hat{u})(\eta)$. Specifically, it is to prove the following theorem holds.
\begin{theorem}\label{th1}
	For $1<\gamma<\frac{3}{2}$, there exists a constant $\bar{\rho}_0>0$, such that for any given initial density $\rho_0\in(0,\bar{\rho}_0)$, there exists a constant $\varepsilon_1>0$, if for any given $\varepsilon\in(0,\varepsilon_1)$,  the piston curve $t=b(r)$ satisfies the following assumptions:
	\begin{itemize}
		\item[\textbf{(H1):}] $b(r) = \frac{r}{b_0}$ for $r \in (0, r_0]$, 
		\item[\textbf{(H2):}] $\left|b'(r) - \frac{1}{b_0}\right| < \frac{\varepsilon}{r^{\delta_0}}$, and $|b''(r)| < \frac{\varepsilon}{r^{1+\delta_0}}$ for $r \in [r_0, +\infty)$,
	\end{itemize}
	for constants $r_0>0$ and $0<\delta_0<3-2\gamma$. Then the piston problem (\ref{e1})-(\ref{e14}) admits a unique piecewise smooth solution with a single shock on the domain $\Omega=\{(r,t)|r> 0,\ t\leqslant b(r)\}$, which has the following form:
	\begin{equation}
		(\rho,u)(r, t) = \begin{cases}
			(\rho_l,u_l)(r, t),\quad & (r,t)\in \Omega_l=\{(r,t)|r>0,\ s(r)\leqslant t\leqslant b(r)\}, \\
			(\rho_0,0), \quad& (r,t)\in \Omega_r=\{(r,t)|r>0,\ t< s(r)\}.
		\end{cases}
	\end{equation}
	Moreover, the following estimates holds on $\Omega_l$:
	\begin{equation}\label{e18}
		\left|(\rho_l- \hat{\rho},u_l - \hat{u})(r,t)\right| \leqslant \mathcal{C} \frac{\varepsilon}{r^{\delta_0}}, \quad \forall (r,t)\in\Omega_l,
	\end{equation}
	and
	\begin{equation}\label{e19}
		\left|(\partial_t(\rho_l- \hat{\rho}),\partial_t(u_l - \hat{u}),\partial_r(\rho_l- \hat{\rho}),\partial_r(u_l - \hat{u}))(r,t)\right| \leqslant \mathcal{C}\frac{\varepsilon}{r^{1+\delta_0}},\quad \forall (r,t)\in\Omega_l.
	\end{equation}
	In addition, the forward shock curve $s(r)\in C^2(\mathbb{R}^+)$ satisfies
	\begin{equation}\label{e10}
		\left|s'(r) - s_0\right|\leqslant \mathcal{C}\frac{\varepsilon}{r^{\delta_0}},\quad\forall r\in\left(0,+\infty\right), 
	\end{equation}
	and
	\begin{equation}\label{e11}
		\left|s''(r)\right|\leqslant \mathcal{C}\frac{\varepsilon}{r^{1+\delta_0}},\quad\forall r\in\left(0,+\infty\right),
	\end{equation}
	where $\mathcal{C}$ is a positive constant independent of $\varepsilon$.
\end{theorem}
\begin{remark}
	The strict dissipation property of leading shock wave degenerates as $\rho_0\to0^+$, see Lemma \ref{le31} for details. Meanwhile, the condition $\gamma<\frac{3}{2}$ is required to balance the shock dissipation and the growth induced by the geometric source term so as to close the a priori estimates in Lemma \ref{le51}. This result corresponds to air with $\gamma=1.4$.
\end{remark}

	The remainder of the paper is organized as follows. In Section 2, we establish some properties for $(\hat{\rho},\hat{u})(\eta)$ and derive relevant geometric estimates for the eigenvalues. In Section 3, we reformulate the original problem as a free boundary problem for perturbation variables via Riemann invariants and then derive the characteristic equations and corresponding boundary conditions for the first order derivatives of perturbation variables. In Section 4, we establish uniform $C^1$ estimates using the reflected characteristics method, following \cite{22,23,24,25,26}, and extend the local solution to a global one via standard continuation arguments \cite{23,24,25} for free boundary problems.

	\section{Preliminary}\label{s2}
	\subsection{Properties of self-similar solutions}
	Before proceeding with the construction for properties of self-simialr solutions, we first clarify the asymptotic notations $O(\cdot)$ and $o(\cdot)$ used throughout this subsection. 
	
    Whenever an infinitesimal quantity $h$ related to $\rho_0$ appears, we uniformly use $O(h)$ to represent any infinitesimal of the same order as $h$, which satisfies
	\begin{equation*}
		C^{-1}h\leqslant O(h) \leqslant C h\quad \text{as } \rho_0 \to 0^+,
	\end{equation*}
	where $C\neq0$ is a positive constant independent of $\rho_0$. Similarly, higher-order infinitesimals of $h$ are denoted by $o(h)$, if
	\begin{equation*}
		\lim\limits_{\rho_0 \to 0^+} \frac{o(h)}{h} = 0.
	\end{equation*}
	
	\begin{lemma}
		For any $1<\gamma<\frac{3}{2}$ and $b_0>0$, we have
		\begin{equation}\label{e29}
			\lim_{\rho_0\to 0^+} s_0 = \frac{1}{b_0}.
		\end{equation}
	Moreover, there holds
	\begin{equation}\label{e216}
		s_0
		=
		\frac1{b_0}
		-
		\frac{1}{3b_0^{1+\frac{2}{\gamma}}}\rho_0^\frac{\gamma-1}{\gamma}
		+
		o(\rho_0^\frac{\gamma-1}{\gamma}),
	\end{equation}
	\end{lemma}
	\begin{proof}
		Introduce the new variables
		\begin{equation*}
			\xi = \eta \hat{u},\qquad \zeta = \eta^2\hat{\rho}^{\gamma-1}.
		\end{equation*}
		On the shock front $\eta=s_0$, denote $\xi_s = s_0 \hat{u}(s_0)$ and $\zeta_s = s_0^2\hat{\rho}(s_0)^{\gamma-1}$, it follows from $(\ref{e21})_4$ and $(\ref{e21})_5$ 
		\begin{equation}\label{e25}
			\zeta_s = \frac{\xi_s(1-\xi_s)}{1-(1-\xi_s)^\gamma},\qquad
			\rho_0 = \Bigl(\frac{\zeta_s}{s_0^2}\Bigr)^{\frac{1}{\gamma-1}}(1-\xi_s).
		\end{equation}
		We can get $\xi_s\to1$ and $\zeta_s\to0$ from (\ref{e25}) as $\rho_0\to0^+$.
		
		A straightforward computation shows
		\begin{equation}\label{e24}
			\frac{d\eta}{\eta} =\frac{d\xi}{\xi}\cdot\frac{(1-\xi)^2-\zeta}{(1-\xi)^2-3\zeta}.
		\end{equation}
		 Integrating the relation (\ref{e24}) from $\left(\frac{1}{b_0},1\right)$ to $\left(s_0,\xi_s\right)$, we obtain
		\begin{equation*}\label{e28}
			\ln s_0=\ln\frac{1}{b_0}+ \int_{1}^{\xi_s}\frac{(1-\hat{\xi})^2-\zeta}{(1-\hat{\xi})^2-3\zeta}\cdot\frac{d\hat{\xi}}{\hat{\xi}}.
		\end{equation*}
		Then as $\rho_0\to0^+$, we get (\ref{e29}). 
		
	Denote $u_s=\hat{u}(s_0),\ \rho_s=\hat{\rho}(s_0)$, it yields
		\begin{equation}\label{eq24}
			u_s=b_0-2b_0^2\left(\frac{1}{b_0}-s_0\right)+o\left(\frac{1}{b_0}-s_0\right),
		\end{equation}
		by Taylor's formula, thus
		\begin{equation}\label{e218}
			u_ss_0=u_s\left[\frac{1}{b_0}-\left(\frac{1}{b_0}-s_0\right)\right]=1-3b_0\left(\frac{1}{b_0}-s_0\right)+o\left(\frac{1}{b_0}-s_0\right).
		\end{equation}
		Substituting (\ref{e218}) into $(\ref{e21})_4$, it yields
		\begin{equation}\label{e219}
			\rho_s=\frac{1}{3b_0}\cdot\frac{\rho_0}{\frac{1}{b_0}-s_0}+o\left(\frac{\rho_0}{\frac{1}{b_0}-s_0}\right).
		\end{equation}
		We can also get the expansion of $\nu=\frac{1}{s_0}$
		\begin{equation*}
			\nu=b_0\left[1+b_0\left(\frac{1}{b_0}-s_0\right)+o\left(\frac{1}{b_0}-s_0\right)\right],
		\end{equation*}
		then
		\begin{equation}\label{eq28}
			\rho_su_s(u_s-\nu)=-\rho_0b_0^2+o(\rho_0).
		\end{equation}
		Since $1<\gamma<\frac{3}{2}$, we have $\rho_0^\gamma=o(\rho_0)$, $(\ref{e21})_5$  and \eqref{eq28} show 
		\begin{equation}\label{e213}
			\rho_s^\gamma=\rho_0b_0^2+o(\rho_0),
		\end{equation}
		 it follows from \eqref{e219} that
		\begin{equation}\label{e222}
			\frac{1}{b_0}-s_0=\frac{1}{3b_0^{1+\frac{2}{\gamma}}}\rho_0^\frac{\gamma-1}{\gamma}+o(\rho_0^\frac{\gamma-1}{\gamma}).
		\end{equation}
		\end{proof}
	
    Denote the sound speed behind the shock $\hat{c}(s_0)$ as $c_s$, \eqref{eq24}, \eqref{e219} and \eqref{e222} show
	\begin{corollary}
		The density, velocity and sound speed behind the shock satisfy
		\begin{align}
			&\rho_s = b_0^{\frac{2}{\gamma}} \rho_0^{\frac{1}{\gamma}} + o(\rho_0^{\frac{1}{\gamma}}),\label{eq323}\\
			&u_s=b_0-\frac{2}{3b_0^{\frac{2}{\gamma}-1}}\rho_0^\frac{\gamma-1}{\gamma}+o(\rho_0^\frac{\gamma-1}{\gamma}),\label{e210}\\
			&c_s =\sqrt{\gamma}\rho_s^\frac{\gamma-1}{2}= \sqrt{\gamma} b_0^{\frac{\gamma-1}{\gamma}}\rho_0^\frac{\gamma-1}{2\gamma} + o(\rho_0^\frac{\gamma-1}{2\gamma})\label{e318}.
		\end{align}
	\end{corollary}

	Next, we establish the monotonicity of self-similar solutions.
	\begin{lemma}\label{l2}
		For any $1<\gamma<\frac{3}{2}$, if $\rho_0$ is small enough, both $\hat{\rho}$ and $\hat{u}$ are monotonically increasing on $\left[s_0,1/b_0\right]$.
	\end{lemma}
	\begin{proof}
		From (\ref{e21}), we have
		\begin{align}
			&\hat{u}_\eta=\frac{2\xi}{\eta^2}\cdot{\frac{\xi}{1-\frac{(1-\xi)^2}{\zeta}}},\label{eq29}\\
			&\hat{\rho}_{\eta}=\frac{1-\xi}{\gamma\eta\rho^{\gamma-2}}\hat{u}_{\eta},\label{eq210}
		\end{align}
		By Taylor's formula and \eqref{e222}, it yields
		\begin{equation}
			\begin{aligned}
				\xi=&\eta\left[b_0-2b_0^2(\eta-s_0)+o(\eta-s_0)\right]\\
				\geqslant&s_0\left[b_0-2b_0^2\left(\frac{1}{b_0}-s_0\right)+o\left(\frac{1}{b_0}-s_0\right)\right]\\
				=&1-b_0^{-\frac{2}{\gamma}}\rho_0^\frac{\gamma-1}{\gamma}+o(\rho_0^\frac{\gamma-1}{\gamma}),
			\end{aligned}
		\end{equation}
		and
		\begin{equation*}
				\zeta=\zeta\left(\frac{1}{b_0}\right)+o(1).
		\end{equation*}
		Therefore, 
		\begin{equation}\label{e 28}
			0\leqslant\frac{(1-\xi)^2}{\zeta}\leqslant\frac{b_0^{-\frac{4}{\gamma}}\rho_0^\frac{2(\gamma-1)}{\gamma}+o(\rho_0^\frac{2(\gamma-1)}{\gamma})}{\zeta\left(\frac{1}{b_0}\right)+o(1)}\to 0,\quad \text{as}\ \rho_0\to0^+.
		\end{equation}
		Furthermore, it follows from \eqref{eq29}-\eqref{e 28} that $\hat{\rho},\ \hat{u}$ are monotically increasing on $\left[s_0,\frac{1}{b_0}\right]$.
	\end{proof}

        Based on the above lemmas, we obtain
    \begin{corollary}
    	 For any $\eta\in\left[s_0,\frac{1}{b_0}\right]$, there holds
    	\begin{align}
    		&b_0-\frac{2}{3b_0^{\frac{2}{\gamma}-1}}\rho_0^\frac{\gamma-1}{\gamma}+o(\rho_0^\frac{\gamma-1}{\gamma})\leqslant\hat{u}(\eta)\leqslant b_0,\label{eq212}\\ 
    		&\hat{u}_{\eta}(\eta)=2b_0^2+O(\rho_0^\frac{\gamma-1}{\gamma})+o(\rho_0^\frac{\gamma-1}{\gamma}),\label{eq213}
    	\end{align}
    	and
    	\begin{align}
    		&\hat{c}(\eta)=\sqrt{\gamma}\hat{\rho}^\frac{\gamma-1}{2}=O(\rho_0^\frac{\gamma-1}{2\gamma})+o(\rho_0^\frac{\gamma-1}{\gamma}).\label{eq214}\\ 
    		&\hat{c}_{\eta}(\eta)=\frac{\sqrt{\gamma}(\gamma-1)}{2}\hat{\rho}^\frac{\gamma-3}{2}\hat{\rho}_{\eta}=O(\rho_0^\frac{\gamma-1}{2\gamma})+o(\rho_0^\frac{\gamma-1}{\gamma}).\label{eq215}
    	\end{align}
    \end{corollary}

	\subsection{Estimates for characteristic curves}
	For application of the reflected characteristics method, a detailed discussion of characteristic curves is required. A direct computation shows that (\ref{e1}) possesses two distinct eigenvalues,
	\begin{equation*}
		\lambda_1(u,c) = u - c,\qquad \lambda_2(u,c) = u + c.
	\end{equation*}
	 We deﬁne the $C_i\ (i = 1,\ 2)$ characteristic curves as the integral curves $\frac{dr}{dt}=\lambda_i$. 
	 \begin{figure}[H]
	 	\centering
	 	\begin{subfigure}[b]{0.48\textwidth}
	 		\centering
	 		\includegraphics[width=\textwidth]{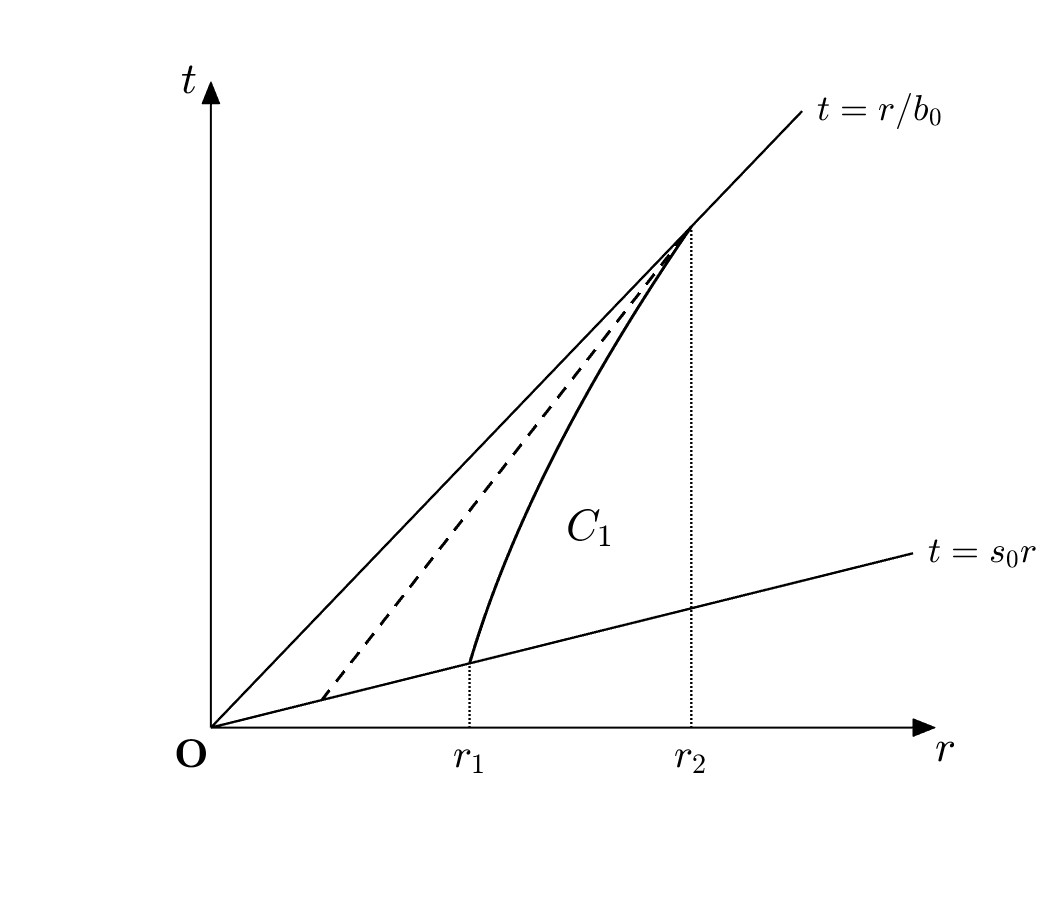}
	 	\end{subfigure}
	 	\hfill 
	 	\begin{subfigure}[b]{0.48\textwidth}
	 		\centering
	 		\includegraphics[width=\textwidth]{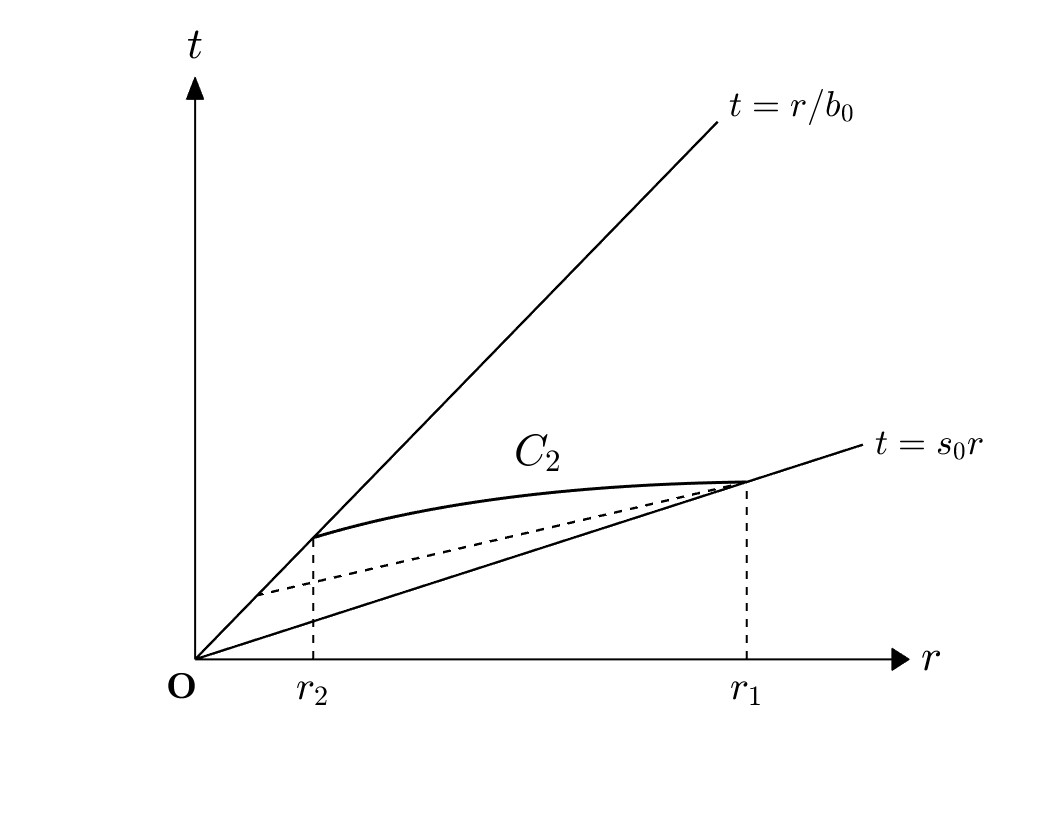}
	 	\end{subfigure}
	 	\caption{$C_1$ and $C_2$ characteristic curves.}
	 	\label{Fig1}
	 \end{figure}
	
	Regarding the estimate for the $C_1$ characteristic, we have
	\begin{lemma}\label{le4}
		Let the $C_1$ characteristic curve emanating from the point $\left(r_2, \frac{r_2}{b_0}\right)$ on the piston intersects the shock front $t=s_0r$ at $(r_1, s_0r_1)$. Then for any given small $\epsilon>0$, when $\rho_0$ is sufficiently small, the following estimate holds
		\begin{equation}\label{e214}
			1 \geqslant \frac{r_1}{r_2} \geqslant 1 - \frac{1}{(3-2\epsilon)\sqrt{\gamma}b_0^{\frac{1}{\gamma}}}\rho_0^\frac{\gamma-1}{2\gamma}.
		\end{equation}
	\end{lemma}
	
	\begin{proof}
		For the $C_1$ characteristic curve, it follows from \eqref{eq212} and \eqref{eq214} that $\lambda_1(\hat{u},\hat{c})>0$, then we have $r_1 \leqslant r_2$, which yields 
		\begin{equation*}
			\frac{r_1}{r_2} \leqslant 1.
		\end{equation*}
		
		Consider the straight line
		\begin{equation*}
			r-r_2=\lambda_1(b_0,c_s)\left(t-\frac{r_2}{b_0}\right).
		\end{equation*}
		By a direct geometric computation, the line emanating from the point $(r_2, r_2/b_0)$ on the piston intersects the shock boundary $t = s_0 r$ at $(\bar{r}_1,s_0\bar{r}_1)$, which is expressed as:
		\begin{equation*}
			\left( r_2 \frac{1 - \frac{\lambda_1(b_0, c_s)}{b_0}}{1 - \lambda_1(b_0, c_s) s_0}, s_0 r_2 \frac{1 - \frac{\lambda_1(b_0, c_s)}{b_0}}{1 - \lambda_1(b_0, c_s) s_0} \right),
		\end{equation*}
		which means
		\begin{equation*}\label{e227}
			\frac{\bar{r}_1}{r_2}=\frac{1 - \frac{\lambda_1(b_0, c_s)}{b_0}}{1 - \lambda_1(b_0,c_s)s_0}.                                    
		\end{equation*}
		Relying on the monotonic properties of $\hat{c}$ and $\hat{u}$, we have
		\begin{equation*}
			\lambda_1(u,c)\leqslant\lambda_1(b_0, c_s),
		\end{equation*}
		thus the true intersection point $(r_1, s_0r_1)$ of the $C_1$ characteristic curve satisfies
		\begin{equation*}
			r_1\geqslant\bar{r}_1.
		\end{equation*}
		Therefore,
		\begin{equation}\label{e228}
			\frac{r_1}{r_2} \geqslant \frac{\bar{r}_1}{r_2}=\frac{1 - \frac{\lambda_1(b_0, c_s)}{b_0}}{1 - \lambda_1(b_0, c_s)s_0}=\frac{1}{1+\frac{\lambda_1(b_0,c_s)(1-b_0s_0)}{b_0-\lambda_1(b_0,c_s)}}.
		\end{equation}
		From \eqref{e216} and \eqref{e318}, we can obtain that for any small $\epsilon>0$, if $\rho_0$ is sufficiently small, then 
		\begin{align*}\label{e230}
			\frac{\lambda_1(b_0,c_s)(1-b_0s_0)}{b_0-\lambda_1(b_0,c_s)}\leqslant \frac{1}{(3-\epsilon)\sqrt{\gamma}b_0^{\frac{1}{\gamma}}}\rho_0^\frac{\gamma-1}{2\gamma}.
		\end{align*}
		Thus, it follows from \eqref{e228}
		\begin{align*}
			\frac{r_1}{r_2} \geqslant\frac{1}{1+\frac{1}{(3-\epsilon)\sqrt{\gamma}b_0^{\frac{1}{\gamma}}}\rho_0^\frac{\gamma-1}{2\gamma}}
			\geqslant1-\frac{1}{(3-2\epsilon)\sqrt{\gamma}b_0^{\frac{1}{\gamma}}}\rho_0^\frac{\gamma-1}{2\gamma}.
		\end{align*}
	\end{proof}
	
	A similar conclusion also holds for $C_2$ characteristic:
	\begin{lemma}\label{le3}
		Assume that the $C_2$ characteristic curve passing through the point $(r_1, s_0r_1)$ on the shock front intersects the piston line $r=b_0t$ at $\left(r_2, \frac{r_2}{b_0}\right)$. Then for any given small $\epsilon>0$, when $\rho_0$ is sufficiently small, the following estimate holds
		\begin{equation}\label{eq 213}
			1 \geqslant \frac{r_2}{r_1} \geqslant 1 - \frac{1}{(3-2\epsilon)\sqrt{\gamma}b_0^{\frac{1}{\gamma}}}\rho_0^\frac{\gamma-1}{2\gamma}.
		\end{equation}
	\end{lemma}
	
	\begin{proof}
	 By the Lax entropy condition, the $C_2$ characteristic speed at the shock front satisfies:
		\begin{equation*}
			\lambda_2(u_s, c_s) = u_s + c_s > \frac{1}{s_0}.
		\end{equation*}
	Then, it yields $r_2 \leqslant r_1$, which immediately gives 
	\begin{equation*}
		\frac{r_2}{r_1} \leqslant 1.
	\end{equation*}
	
	The $C_2$ characteristic curve in the $t$-$r$ plane satisfies the differential equation:
	\begin{equation}\label{eq219}
		\frac{dt}{dr} = \frac{1}{\lambda_2(\hat{u}, \hat{c})} = \frac{1}{\hat{u} + \hat{c}}.
	\end{equation}
	It follows from Lemma \ref{l2} that
	\begin{equation*}
		\frac{1}{\lambda_2(\hat{u}, \hat{c})} \leqslant \frac{1}{\lambda_2(u_s, c_s)}.
	\end{equation*}
	Then integrating (\ref{eq219}) from $r_2$ to $r_1$ yields
		\begin{equation*}
			s_0 r_1 - \frac{r_2}{b_0} = \int_{r_2}^{r_1} \frac{1}{\lambda_2(\hat{u}, \hat{c})} dr \leqslant \int_{r_2}^{r_1} \frac{1}{\lambda_2(u_s, c_s)} dr = \frac{r_1 - r_2}{\lambda_2(u_s, c_s)},
		\end{equation*}
		thus
		\begin{equation*}\label{eq222}
			\frac{r_2}{r_1} \geqslant \frac{\frac{1}{\lambda_2(u_s, c_s)} - s_0}{\frac{1}{\lambda_2(u_s, c_s)} - \frac{1}{b_0}} = b_0 s_0 \frac{\lambda_2(u_s, c_s) s_0 - 1}{s_0 (\lambda_2(u_s, c_s) - b_0)} = b_0 s_0 \left(1 - \frac{\frac{1}{s_0} - b_0}{\lambda_2(u_s, c_s) - b_0}\right).
		\end{equation*}
		By \eqref{e216}, \eqref{e210} and \eqref{e318} we obtain that for any small $\epsilon>0$, if $\rho_0$ is sufficiently small, then 
		\begin{equation*}
			\frac{r_2}{r_1} \geqslant \left(1 - \frac{1}{(3-\epsilon)b_0^{\frac{2}{\gamma}}}\rho_0^\frac{\gamma-1}{\gamma}\right)\cdot \left( 1 - \frac{1}{(3-\epsilon)\sqrt{\gamma}b_0^{\frac{1}{\gamma}}}\rho_0^\frac{\gamma-1}{2\gamma} \right),
		\end{equation*}
		which shows 
		\begin{equation*}
			\frac{r_2}{r_1} \geqslant 1 - \frac{1}{(3-2\epsilon)\sqrt{\gamma}b_0^{\frac{1}{\gamma}}}\rho_0^\frac{\gamma-1}{2\gamma}.
		\end{equation*}   
		Consequently, the proof is completed.
	\end{proof}

	In the region $\Omega=\left\{(r,t)\left\lvert r>0,\ s_0r\leqslant t\leqslant\frac{r}{b_0}\right.\right\}$, we have the following pointwise estimates for the characteristic speeds.
	\begin{lemma}\label{le5}
		 When $\rho_0>0$ is sufficiently small, the following estimate holds for $j=1,\ 2$
	 \begin{equation}\label{e233}
	 |b_0 - \lambda_j(\hat{u}, \hat{c})| =  O(\rho_0^\frac{\gamma-1}{2\gamma}).
	 \end{equation}
    \end{lemma}
    \begin{proof}
      For $j=1$, it follows from \eqref{eq212} and \eqref{eq214}
      \begin{equation*}\label{eq227}
    	|b_0 - \lambda_1(\hat{u}, \hat{c})| = |b_0 - \hat{u} + \hat{c}| = (b_0 - \hat{u}) + \hat{c} \geqslant \hat{c}(\eta) = O(\rho_0^\frac{\gamma-1}{2\gamma})+o(\rho_0^\frac{\gamma-1}{\gamma}).
      \end{equation*}
      For $j=2$, 
      \begin{equation*}
      	|b_0 - \lambda_2(\hat{u}, \hat{c})| = |b_0 - \hat{u} - \hat{c}| = |\hat{c} - (b_0 - \hat{u})|.
      \end{equation*}
      Since $b_0 - \hat{u} = O(\rho_0^\frac{\gamma-1}{\gamma})$ and $\hat{c} = O(\rho_0^\frac{\gamma-1}{2\gamma})+o(\rho_0^\frac{\gamma-1}{\gamma})$, $O(\rho_0^\frac{\gamma-1}{\gamma})$ can be absorbed by $O(\rho_0^\frac{\gamma-1}{2\gamma})$ for sufficiently small $\rho_0$. Therefore, we have
      \begin{equation*}\label{eq228}
      	|b_0 - \lambda_2(\hat{u}, \hat{c})| =  O(\rho_0^\frac{\gamma-1}{2\gamma}).
      \end{equation*}

    \end{proof}
    
    \section{Reformulation of the problem}\label{s3}
    We aim to construct a global piecewise smooth solution containing a single shock. To this end, we regard the solution of  \eqref{e1}-\eqref{e14} as a small perturbation of the background self-similar solutions $(\hat{\rho}, \hat{u})$ obtained in Section 2.
    In this section, we will derive the diagonalized equations and boundary conditions satisfied by the perturbation variables.
    
    \subsection{The derivation of equivalent free boundary problem}
    Consider system (\ref{e1}) with the fixed boundary condition on the piston
    \begin{equation}\label{eq31}
    	u_l(r, t) = \frac{1}{b'(r)} \quad \text{on } \Gamma_b,
    \end{equation}
    where $\Gamma_b := \{(r, t) \mid r \geqslant r_0, t = b(r)\}$, and the free boundary condition on the shock
    \begin{equation}\label{eq32}
    	H(\rho_l,u_l) = 0 \quad \text{and} \quad s'(r) = \bar{\nu}(u_l) \quad \text{on } \Gamma_s,
    \end{equation}
    where $\Gamma_s := \{(r, t) \mid r \geqslant r_0, t = s(r)\}$. And from the Rankine-Hugoniot conditions (\ref{e14}), it yields
    \begin{equation}\label{e33}
    	H(\rho_l,u_l) := \rho_0 \rho_l u_l^2 - (\rho_l - \rho_0)(\rho_l^\gamma - \rho_0^\gamma) = 0,
    \end{equation}
     and
    \begin{equation}\label{eq33}
    	\bar{\nu}(u_l) = \frac{\rho_l - \rho_0}{\rho_l u_l}.
    \end{equation}
    
    By directly calculating for \eqref{eq33}, we obtain the following lemma.
    \begin{lemma}\label{le 31}
    	 $\bar{\nu}'(u_l)$ satisfies the following conditions:
    	\begin{align}
    		\lim\limits_{\rho_0\to 0^+}\bar{\nu}'(u_l)=-\frac{1}{u_l^2}.\label{e34}
    	\end{align}
    \end{lemma}

    Introduce the Riemann invariants
    \begin{equation}\label{n36}
    	m = u_l - \frac{2c_l}{\gamma - 1}, \quad n = u_l + \frac{2c_l}{\gamma - 1}.
    \end{equation}
    Then system (\ref{e1}) can be diagonalized as
    \begin{equation}\label{eq34}
    	\begin{cases}
    		\partial_t m + \lambda_1(m, n) \partial_r m = \frac{F(m,n)}{r}, \\
    		\partial_t n + \lambda_2(m, n) \partial_r n = -\frac{F(m,n)}{r},
    	\end{cases}
    \end{equation}
    with
    \begin{align*}
    	&\lambda_1(m,n) = \frac{1+\gamma}{4}m+\frac{3-\gamma}{4}n,\\ &\lambda_2(m,n)= \frac{3-\gamma}{4}m+\frac{1+\gamma}{4}n,\\
    	 &F(m,n) = \frac{\gamma-1}{4}\left(n^2-m^2\right).
    \end{align*}
    
    Define the background Riemann invariants based on the self-similar solution $(\hat{u}, \hat{c})(\eta)$ with $\eta = \frac{t}{r}$
    \begin{equation}\label{n38}
    	\hat{m}(r, t) = \hat{u}(\eta) - \frac{2\hat{c}(\eta)}{\gamma - 1}, \quad \hat{n}(r, t) = \hat{u}(\eta) + \frac{2\hat{c}(\eta)}{\gamma - 1},
    \end{equation}
    then $\eqref{e21}_1-\eqref{e21}_2$ become
    	\begin{equation}\label{eq35}
    		\begin{cases}
    			\partial_t \hat{m} + \lambda_1(\hat{m}, \hat{n}) \partial_r \hat{m} = \frac{F(\hat{m},\hat{n})}{r}, \\
    			\partial_t \hat{n} + \lambda_2(\hat{m}, \hat{n}) \partial_r \hat{n} = -\frac{F(\hat{m},\hat{n})}{r}.
    		\end{cases}
    	\end{equation}
    Furthermore, we have the following estimates.
     \begin{lemma}\label{le30}
    	The variables $\partial_r\hat{m}$ and $\partial_r\hat{n}$ satisfy the following conditions:
    	\begin{equation}\label{e313}
    		\lim_{\rho_0\to 0^+}\left|\partial_r\hat{m}\right|=\lim_{\rho_0\to 0^+}\left|\partial_r\hat{n}\right|=\frac{2}{r}b_0.
    	\end{equation}
    \end{lemma}
    \begin{proof}
    	By computation, we have
    	\begin{equation*}
    		\partial_r\hat{m}=-\frac{\eta}{r}\left(\hat{u}_{\eta}-\frac{2}{\gamma-1}\hat{c}_{\eta}\right),
    	\end{equation*}
    	and
    	\begin{equation*}
    		\partial_r\hat{n}=-\frac{\eta}{r}\left(\hat{u}_{\eta}+\frac{2}{\gamma-1}\hat{c}_{\eta}\right).
    	\end{equation*}
    	Therefore, using \eqref{eq213} and \eqref{eq215}, we get that \eqref{e313} holds.
    \end{proof}
    
    Introduce the perturbation variables
    \begin{equation}\label{eq36}
    	\bar{m}(r, t) = m(r, t) - \hat{m}(r, t), \quad \bar{n}(r, t) = n(r, t) - \hat{n}(r, t),
    \end{equation}
    then, by \eqref{eq34}-\eqref{eq35}, we obtain
    \begin{equation}\label{eq37}
    	\begin{cases}
    		\partial_t \bar{m} + \lambda_1(m, n) \partial_r \bar{m} = \frac{F(m, n) - F(\hat{m}, \hat{n})}{r} + \big[\lambda_1(\hat{m}, \hat{n}) - \lambda_1(m, n)\big]\partial_r \hat{m}, \\
    		\partial_t \bar{n} + \lambda_2(m, n) \partial_r \bar{n} = -\frac{F(m, n) - F(\hat{m}, \hat{n})}{r} + \big[\lambda_2(\hat{m}, \hat{n}) - \lambda_2(m, n)\big]\partial_r \hat{n},
    	\end{cases}
    \end{equation}
    where
    \begin{align*}\label{eq38}
    		&F(m, n) - F(\hat{m}, \hat{n})= \frac{\gamma-1}{2}(\hat{n}\bar{n}-\hat{m}\bar{m}) + \frac{\gamma-1}{4}(\bar{n}^2-\bar{m}^2), \\
    		&\lambda_1(m, n) - \lambda_1(\hat{m}, \hat{n})= \frac{\gamma+1}{4} \bar{m} + \frac{3-\gamma}{4} \bar{n} ,\\
    	&\lambda_2(m, n) - \lambda_2(\hat{m}, \hat{n})= \frac{3-\gamma}{4} \bar{m} + \frac{\gamma+1}{4} \bar{n} .
    \end{align*}

    Next, we give the boundary conditions of $\bar{m}$ and $\bar{n}$. First, since $\rho$ and $u$ can be expressed by $m, n$, the first equality of \eqref{eq32} can be written as
    \begin{equation}\label{eq39}
    	H(\bar{m} + \hat{m},\bar{n} + \hat{n}) = 0 \quad \text{on } \Gamma_s.
    \end{equation}
    Moreover, since $\hat{m},\hat{n}$ are functions of $\frac{t}{r}$, we obtain the boundary condition on the shock curve $\Gamma_s$ is
    \begin{equation}\label{eq310}
    	\hat{H}(\bar{m},\bar{n},\frac{t}{r}):=H(\bar{m} + \hat{m},\bar{n} + \hat{n})=0.
    \end{equation}
    Define 
    \begin{equation*}\label{eq:3.15}
    	\tau=\tau(r) = \frac{s(r)}{r} - s_0, \quad r > 0.
    \end{equation*}
    Since $\hat{H}(0,0,s_0) = 0$, applying Taylor's formula and implicit function theorem, the boundary condition \eqref{eq310} can be linearized as:
    \begin{equation}\label{eq316}
    	\bar{m}(r, t) = A_1 \bar{n}(r, t) + A_2 \tau(r) + I(\bar{n}(r, t), \tau(r)) \quad \text{on } \Gamma_s,
    \end{equation} 
    where
    \begin{equation}\label{eq:3.18}
    	A_1 = -\frac{l_2}{l_1}, \quad A_2 = -\frac{\hat{m}'(s_0)l_1 + \hat{n}'(s_0)l_2}{l_1},
    \end{equation}
    and
    \begin{equation*}
    	l_1 = \rho_0 \rho_s u_s - \frac{\rho_s}{2c_s} \Big[ \rho_0 u_s^2 - (\rho_s^\gamma - \rho_0^\gamma) - c^2(\rho_s - \rho_0) \Big],
    \end{equation*}
    \begin{equation*}
    	l_2 = \rho_0 \rho_s u_s + \frac{\rho_s}{2c_s} \Big[ \rho_0 u_s^2 - (\rho_s^\gamma - \rho_0^\gamma) - c_s^2(\rho_s - \rho_0) \Big].
    \end{equation*}
     The non-linear remainder $I(\bar{n}, \tau)$ satisfies $I(0,0) = \partial_{\bar{n}}I(0,0) = \partial_\tau I(0,0) = 0$, then there exists a positive constant $\mathcal{M}$ independent of $\bar{n}$ and $\tau$, such that as $\bar{n}$ and $\tau$ are sufficiently small,
     \begin{equation}\label{eq314}
     	\left\{
     	\begin{aligned}
     		&\left\lvert I(\bar{n},\tau)\right\rvert\leqslant \mathcal{M}\left(|\bar{n}|^2+|\tau|^2\right),\\
     	&\left\lvert \partial_{\bar{n}}I(\bar{n},\tau)\right\rvert\leqslant \mathcal{M}\left(|\bar{n}|+|\tau|\right),\\
     	&\left\lvert \partial_{\tau}I(\bar{n},\tau)\right\rvert\leqslant \mathcal{M}\left(|\bar{n}|+|\tau|\right).
     	\end{aligned}
     	\right.
     \end{equation}
    
    By the second equality of \eqref{eq32}, we have
    \begin{equation}\label{eq317}
    	(r \tau(r))' = \bar{\nu}(u(r, s(r))) - s_0,
    \end{equation}
    then exist $\theta_k\in\left(0,1\right)\ \left(k=1,2,3\right)$, such that (\ref{eq317}) can be written as
    \begin{align}\label{e320}
    	(r \tau(r))' =& \bar{\nu}(u_l(r, s(r))) - s_0,\notag\\
    	=&\bar{\nu}(u_l(r, s(r)))-\bar{\nu}\left(\hat{u}(\frac{s(r)}{r})\right)+\bar{\nu}\left(\hat{u}(\frac{s(r)}{r})\right)-\bar{\nu}\left(\hat{u}(s_0)\right)\notag\\
    	=&\bar{\nu}'\left(\hat{u}(\frac{s(r)}{r})+\theta_1\left(u_l(r,s(r))-\hat{u}(\frac{s(r)}{r})\right)\right)\left(u_l(r,s(r))-\hat{u}(\frac{s(r)}{r})\right)\\
    	&+\bar{\nu}'\left(\hat{u}(s_0)+\theta_2\left(\hat{u}(\frac{s(r)}{r})-\hat{u}(s_0)\right)\right)\left(\hat{u}(\frac{s(r)}{r})-\hat{u}(s_0)\right)\notag\\
    	=&\bar{\nu}'\left(\hat{u}(s_0+\tau)+\theta_1\left(u_l(r,s(r))-\hat{u}(\frac{s(r)}{r})\right)\right)\cdot\left(u_l(r,s(r))-\hat{u}(\frac{s(r)}{r})\right)\notag\\
    	&+\bar{\nu}'\left(\hat{u}(s_0)+\theta_2\hat{u}_\eta(s_0+\theta_3\tau)\tau\right)\hat{u}_\eta(s_0+\theta_3\tau)\tau.\notag
    \end{align}

    \begin{lemma} \label{le31}
    	For $A_1$ and $A_2$, we have
    \begin{align}
    	|A_1| &= 1 -  \frac{4}{\sqrt{\gamma}}b_0^{-\frac{1}{\gamma}}  \rho_0^\frac{\gamma-1}{2\gamma}+o(\rho_0^\frac{\gamma-1}{2\gamma}) , \label{eq3161} \\
    	|A_2| &= \frac{4}{\sqrt{\gamma}}b_0^{2-\frac{1}{\gamma}} \rho_0^\frac{\gamma-1}{2\gamma}+o(\rho_0^\frac{\gamma-1}{2\gamma}) . \label{3162}
    \end{align}
    \end{lemma}
    
    \begin{proof}
    	We first prove the estimate for $A_1$. Set variables $X = \rho_0\rho_s u_s$ and $Y= \frac{\rho_s}{2c_s} \Big[ \rho_0 u_s^2 - (\rho_s^\gamma - \rho_0^\gamma) - c_s^2(\rho_s - \rho_0) \Big]$. By \eqref{e33}, we have
    	\begin{equation}\label{e317}
    		-Y = \frac{\rho_s}{2c_s} \left[ \frac{\rho_0}{\rho_s}(\rho_s^\gamma - \rho_0^\gamma) + c_s^2(\rho_s - \rho_0) \right].
    	\end{equation}
    	From \eqref{eq323}-\eqref{e318},
    	\begin{align*}
    		X&=b_0^\frac{\gamma+2}{\gamma}\rho_0^\frac{\gamma+1}{\gamma}+o(\rho_0^\frac{\gamma+1}{\gamma}),\\
    		-Y&=\frac{\rho_s}{2c_s} \left[\frac{\rho_0}{\rho_s}(\rho_s^\gamma - \rho_0^\gamma) + c_s^2(\rho_s - \rho_0) \right]
    		=\frac{\sqrt{\gamma}}{2}b_0^\frac{\gamma+3}{\gamma}\rho_0^\frac{\gamma+3}{2\gamma}+o(\rho_0^\frac{\gamma+3}{2\gamma}),
    	\end{align*}
    	thus 
    	\begin{equation*}
    		\frac{X}{-Y} = \frac{2}{\sqrt{\gamma}}b_0^{-\frac{1}{\gamma}}  \rho_0^\frac{\gamma-1}{2\gamma}+o(\rho_0^\frac{\gamma-1}{2\gamma}).
    	\end{equation*}
    	Hence, we obtain:
    	\begin{equation}\label{e329}
    		A_1 = \frac{1-\frac{X}{-Y}}{1+\frac{X}{-Y}} = 1 -  \frac{4}{\sqrt{\gamma}}b_0^{-\frac{1}{\gamma}}  \rho_0^\frac{\gamma-1}{2\gamma}+o(\rho_0^\frac{\gamma-1}{2\gamma}),
    	\end{equation}
    	which leads to \eqref{eq3161}.
    	
    	Next, we prove the estimate for $A_2$. By definition \eqref{eq:3.18}, $A_2 = -\hat{m}'(s_0) + A_1 \hat{n}'(s_0)$.
    	Using \eqref{e329}, we have:
    	\begin{align*}
    		A_2 &= -\left(\hat{u}_\eta - \frac{2}{\gamma-1}\hat{c}_\eta\right) + \left(1 - \frac{4}{\sqrt{\gamma}}b_0^{-\frac{1}{\gamma}}  \rho_0^\frac{\gamma-1}{2\gamma}+o(\rho_0^\frac{\gamma-1}{2\gamma}) \right)\cdot\left(\hat{u}_\eta + \frac{2}{\gamma-1}\hat{c}_\eta\right) \\
    		&= \frac{4}{\gamma-1}\left(1-\frac{2}{\sqrt{\gamma}}b_0^{-\frac{1}{\gamma}}  \rho_0^\frac{\gamma-1}{2\gamma}+o(\rho_0^\frac{\gamma-1}{2\gamma})\right)\hat{c}_\eta - \left(\frac{4}{\sqrt{\gamma}}b_0^{-\frac{1}{\gamma}}  \rho_0^\frac{\gamma-1}{2\gamma}+o(\rho_0^\frac{\gamma-1}{2\gamma})\right)\hat{u}_\eta.
    	\end{align*}
    	We now evaluate $\hat{c}_\eta(s_0)$. From $(\ref{e21})_2$, we deduce $\frac{\hat{\rho}_\eta}{\hat{\rho}} = \frac{1-s_0 \hat{u}}{s_0 \hat{c}^2}\hat{u}_\eta$, together with $2 \hat{c} \hat{c}_\eta = \gamma(\gamma-1)\hat{\rho}^{\gamma-2}\hat{\rho}_\eta$, we have
    	\begin{equation}\label{e325}
    		\hat{c}_\eta(s_0) = \frac{\gamma-1}{2 c_s s_0} (1 - s_0 u_s) \hat{u}_\eta(s_0).
    	\end{equation}
    	Substituting (\ref{e218}), (\ref{e318}) and (\ref{eq213}) into (\ref{e325}), we obtain
    	\begin{equation*}
    		\hat{c}_\eta(s_0)= \frac{\gamma-1}{\sqrt{\gamma}}b_0^{2-\frac{1}{\gamma}}\rho_0^\frac{\gamma-1}{2\gamma}+o(\rho_0^\frac{\gamma-1}{2\gamma}).
    	\end{equation*}
    	Therefore, we have
    	\begin{equation*}
    		A_2 = -\frac{4}{\sqrt{\gamma}}b_0^{2-\frac{1}{\gamma}}\rho_0^\frac{\gamma-1}{2\gamma}+o(\rho_0^\frac{\gamma-1}{2\gamma}),
    	\end{equation*}
    	which leads to \eqref{3162}. We complete the proof of Lemma \ref{le31}.
    \end{proof}

    Finally, we reformulate the fixed boundary condition (\ref{eq31}) on the piston $\Gamma_b$. Using $u = \frac{m + n}{2}$, we have:
    \begin{equation*}\label{eq:3.19}
    	\bar{m}(r, t) + \bar{n}(r, t) + \hat{m}\Big(\frac{t}{r}\Big) + \hat{n}\Big(\frac{t}{r}\Big) = \frac{2}{b'(r)} \quad \text{on } \Gamma_b,
    \end{equation*}
    which can be explicitly rewritten as
    \begin{equation}\label{eq320}
    	\bar{n}(r, t) = -\bar{m}(r, t) + G_{11}(r) + G_{12}(r, t) \quad \text{on } \Gamma_b,
    \end{equation}
    where 
    \begin{equation*}
    	G_{11}(r) = 2\left(\frac{1}{b'(r)} - b_0\right), \quad G_{12}(r, t) = 2b_0 - \left(\hat{m}\Big(\frac{t}{r}\Big) + \hat{n}\Big(\frac{t}{r}\Big)\right).
    \end{equation*}
    
    \begin{lemma}\label{le32}
    	Under assumptions \textbf{(H1)}--\textbf{(H2)}, when $\varepsilon$ is sufficiently small, $G_{11}(r)$ and $G_{12}(r, t)$ satisfy the following decay estimates:
    	\begin{equation}\label{eq:3.22}
    		\begin{cases}
    			|G_{11}(r)| \leqslant 4b_0^2\cdot\frac{\varepsilon}{r^{\delta_0}}, & r > r_0; \\
    			|G_{12}(r, t)| \leqslant \frac{5b_0^2}{1-\delta_0}\cdot\frac{\varepsilon}{r^{\delta_0}}, & (r, t) \in \Gamma_b.
    		\end{cases}
    	\end{equation}
    	
    \end{lemma}
    \begin{proof}
    	The bound for $G_{11}(r)$ follows 
    	\begin{equation*}
    		|G_{11}(r)|=\left|\frac{2b_0}{b'(r)}\left(b'(r)-\frac{1}{b_0}\right)\right|\leqslant\frac{2b_0}{\frac{1}{b_0}-\frac{\varepsilon}{r^{\delta_0}}}\cdot\frac{\varepsilon}{r^{\delta_0}}\leqslant4b_0^2\cdot\frac{\varepsilon}{r^{\delta_0}}.
    	\end{equation*}
    	
    	 For $G_{12}(r, t)$, noting that $\hat{m}(\frac{1}{b_0}) + \hat{n}(\frac{1}{b_0}) = 2\hat{u}(\frac{1}{b_0}) = 2b_0$, we have
    	\begin{equation}\label{eq336}
    		|G_{12}(r, t)| = \left| \left(\hat{m}\Big(\frac{1}{b_0}\Big) + \hat{n}\Big(\frac{1}{b_0}\Big)\right) - \left(\hat{m}\Big(\frac{t}{r}\Big) + \hat{n}\Big(\frac{t}{r}\Big)\right) \right| \leqslant M \left| \frac{1}{b_0} - \frac{b(r)}{r} \right|,
    	\end{equation}
    	where $M = \sup\limits_{\eta\in\left[s_0,\frac{1}{b_0}\right]} |\hat{m}_\eta(\eta) + \hat{n}_\eta(\eta)|\leqslant5b_0^2$, and 
    	\begin{equation}\label{eq334}
    		\left|\frac{1}{b_0} - \frac{b(r)}{r}\right|=\left\lvert\frac{1}{r}\left(\frac{r}{b_0}-b(r)\right)\right\rvert=\left\lvert\frac{1}{r}\int_{r_0}^{r}\left(\frac{1}{b_0}-b'(z)\right)dz\right\rvert<\frac{1}{1-\delta_0}\cdot\frac{\varepsilon}{r^{\delta_0}}.
    	\end{equation}
        Substituting \eqref{eq334} into \eqref{eq336}, we obtain $\eqref{eq:3.22}_2$.
    \end{proof}

    \subsection{Characteristic equations and boundary conditions for the derivatives}
    
    To establish the uniform $C^1$ estimates for the classical solution, we need to derive the characteristic equations and boundary conditions for the first-order derivatives of the perturbations $\bar{m}$ and $\bar{n}$. 
    
    Following the geometric framework, we introduce a specific directional derivative along the unperturbed piston trajectory:
    \begin{equation*}\label{eq41}
    	\partial_p = \partial_r + \frac{1}{b_0}\partial_t.
    \end{equation*}
    We define the new variables representing the derivatives of the perturbations as:
    \begin{equation*}
    	\dot{m} = \partial_p \bar{m}, \quad \dot{n}=\partial_p \bar{n}.
    \end{equation*}

    Let $\partial_j = \partial_t + \lambda_j \partial_r$ ($j=1, 2$), the partial derivatives $\partial_r$ and $\partial_t$ can be represented by $\partial_p$ and $\partial_j$:
    \begin{equation}\label{eq43}
    	\partial_r = \frac{\partial_j - b_0 \partial_p}{\lambda_j - b_0} \quad \text{and} \quad \partial_t = \frac{\partial_j - \lambda_j \partial_p}{1 - \frac{\lambda_j}{b_0}}.
    \end{equation}
    Applying the operator $\partial_p$ to the perturbed system (\ref{eq37}), we obtain the characteristic equations for $\dot{m}$ and $\dot{n}$:
    \begin{equation}\label{eq44}
    	\begin{cases}
    		\partial_t \dot{m} + \lambda_1(m, n) \partial_r \dot{m} = H_1(m,n), \\
    		\partial_t \dot{n} + \lambda_2(m, n) \partial_r \dot{n} = H_2(m,n),
    	\end{cases}
    \end{equation}
    where 
    \begin{align}
    	H_1(m,n) &= \frac{1}{r} \partial_p \big[ F(m, n) - F(\hat{m}, \hat{n}) \big] + \partial_p \big[ \lambda_1(\hat{m}, \hat{n}) - \lambda_1(m, n) \big](\partial_r \hat{m} + \partial_r \bar{m}) \notag \\
    	&\quad + \big[ \lambda_1(\hat{m}, \hat{n}) - \lambda_1(m, n) \big]\partial_p \partial_r \hat{m} - \partial_p \lambda_1(\hat{m}, \hat{n}) \partial_r \bar{m} - \frac{F(m, n) - F(\hat{m}, \hat{n})}{r^2}, \label{eq45} \\
    	H_2(m,n) &= -\frac{1}{r} \partial_p \big[ F(m, n) - F(\hat{m}, \hat{n}) \big] + \partial_p \big[ \lambda_2(\hat{m}, \hat{n}) - \lambda_2(m, n) \big](\partial_r \hat{n} + \partial_r \bar{n}) \notag \\
    	&\quad + \big[ \lambda_2(\hat{m}, \hat{n}) - \lambda_2(m, n) \big]\partial_p \partial_r \hat{n} - \partial_p \lambda_2(\hat{m}, \hat{n}) \partial_r \bar{n} + \frac{F(m, n) - F(\hat{m}, \hat{n})}{r^2}. \label{eq46}
    \end{align}

    Moreover, we can compute the $\partial_p$ derivatives explicitly
    \begin{align}
    	&\partial_p \big[ \lambda_1(m, n) - \lambda_1(\hat{m}, \hat{n}) \big] = \frac{\gamma+1}{4}\dot{m}+\frac{3-\gamma}{4}\dot{n}, \label{eq47} \\
    	&\partial_p \big[ \lambda_2(m, n) - \lambda_2(\hat{m}, \hat{n}) \big] = \frac{3-\gamma}{4}\dot{m}+\frac{\gamma+1}{4}\dot{n}, \label{eq48}\\
    	&\partial_p \big[ F(m, n) - F(\hat{m}, \hat{n}) \big] = -\frac{\gamma-1}{2}\big( \hat{m} \dot{m} - \hat{n} \dot{n} + \bar{m} \partial_p \hat{m} - \bar{n} \partial_p \hat{n} \big) - \frac{\gamma-1}{2}\big( \bar{m} \dot{m} - \bar{n} \dot{n} \big).
    	\end{align}
    Using (\ref{eq43}), one has
    \begin{equation}\label{eq411}
    	\partial_r \bar{m} = \frac{1}{b_0 - \lambda_1} \left( b_0 \dot{m} - \partial_1 \bar{m} \right),\ \partial_r \bar{n} = \frac{1}{b_0 - \lambda_2} \left( b_0 \dot{n} - \partial_2 \bar{n} \right)
    \end{equation}
    These explicit representations \eqref{eq45}-\eqref{eq411} ensure that $H_i$ are strictly bounded by combinations of $\dot{m}$, $\dot{n}$, $\bar{m}$, $\bar{n}$, $\hat{m}$ and $\hat{n}$, preparing a foundation for the uniform estimates. Besides, a direct calculation shows
    \begin{lemma}
    	For any small constant $\epsilon>0$, when $\rho_0$ is sufficiently small, the following conclusion holds consistently on $[s_0,1/b_0]$:
    	\begin{equation}\label{e337}
    	0<\partial_p\hat{m},\ 	\partial_p\hat{n}\leqslant\frac{1}{r}\cdot\left(\frac{2}{3}+\epsilon\right)b_0^{1-\frac{2}{\gamma}}\rho_0^\frac{\gamma-1}{\gamma}.
    	\end{equation}
    \end{lemma}

    Next, we derive the boundary conditions satisfied by $\dot{m}$ and $\dot{n}$ on $\Gamma_b$. Let $\partial_b = \partial_r + b'(r)\partial_t$ denote the tangential derivative along the perturbed piston curve $\Gamma_b$. Differentiating the fixed boundary condition (\ref{eq320}) along $\Gamma_b$, we have
    \begin{equation}\label{eq412}
    	\partial_b \bar{n} = -\partial_b \bar{m} + G'_{11}(r) + \partial_b G_{12}(r, t) \quad \text{on } t = b(r).
    \end{equation}
    To relate $\partial_b \bar{m}$ and $\partial_b \bar{n}$ to our target variable $\dot{m}$ and $\dot{n}$, we substitute the identity (\ref{eq43}) to get
    \begin{equation}\label{eq413}
    	\partial_b \bar{m} = \frac{1 - b'(r)\lambda_1}{1 - \frac{\lambda_1}{b_0}} \dot{m} + \frac{b'(r) - \frac{1}{b_0}}{1 - \frac{\lambda_1}{b_0}} \partial_1 \bar{m},\quad \partial_b \bar{n} = \frac{1 - b'(r)\lambda_2}{1 - \frac{\lambda_2}{b_0}} \dot{n} + \frac{b'(r) - \frac{1}{b_0}}{1 - \frac{\lambda_2}{b_0}} \partial_2 \bar{n}.
    \end{equation}
    Substituting \eqref{eq413} into \eqref{eq412}, we establish the boundary condition for $\dot{n}$ on $\Gamma_b$:
    \begin{align}\label{eq:4.12}
    	\dot{n} &= -\dot{m} + \left( b'(r) - \frac{1}{b_0} \right) \left( \frac{\lambda_2 \dot{n}}{1 - \frac{\lambda_2}{b_0}} + \frac{\lambda_1 \dot{m}}{1 - \frac{\lambda_1}{b_0}} \right) 
    	 - \frac{b'(r) - \frac{1}{b_0}}{1 - \frac{\lambda_2}{b_0}} \partial_2 \bar{n} - \frac{b'(r) - \frac{1}{b_0}}{1 - \frac{\lambda_1}{b_0}} \partial_1 \bar{m} \notag\\
    	&\quad - \Big[ \frac{2b''(r)}{1-[b'(r)]^2} + \frac{1}{r}\left(b'(r)-\frac{t}{r}\right)\left(\hat{m}_\eta+\hat{n}_\eta\right)\Big] \quad \text{on } \Gamma_b.
    \end{align}
    By assumption \textbf{(H2)}, the terms containing $b'(r) - \frac{1}{b_0}$ in \eqref{eq:4.12} are higher-order smallness terms of $\frac{\varepsilon}{r^{\delta_0}}$. Moreover, we have
    \begin{align}\label{eq416}
    	\left\lvert b'(r)-\frac{t}{r}\right\rvert&\leqslant\left\lvert b'(r)-\frac{1}{b_0}\right\rvert+\left\lvert\frac{1}{b_0}-\frac{t}{r}\right\rvert\notag\\
    	&=\left\lvert b'(r)-\frac{1}{b_0}\right\rvert+\frac{1}{r}\left\lvert\frac{r}{b_0}-b(r)\right\rvert\\
    	&<\left(1+\frac{1}{1-\delta_0}\right)\frac{\varepsilon}{r^{\delta_0}}\notag
    \end{align}
    through (\ref{eq334}).
    
    Finally, we deduce the boundary conditions satisfied by $\dot{m}$ and $\dot{n}$ on $\Gamma_s$. Similarly, let $\partial_s = \partial_r + s'(r)\partial_t$ denote the tangential derivative along the shock front $\Gamma_s$. Differentiating the free boundary condition \eqref{eq316} along $\Gamma_s$ yields:
    \begin{equation}\label{eq:4.13}
    	\partial_s \bar{m} = A_1 \partial_s \bar{n} + A_2 \tau'(r) + \partial_{\bar{n}} I \partial_s \bar{n} + \partial_\tau I \tau'(r) \quad \text{on } \Gamma_s.
    \end{equation}
    Following the same coordinate transformation as (\ref{eq413}) but replacing $b'(r)$ with $s'(r)$, and then substitute it into \eqref{eq:4.13} to obtain
    \begin{align}\label{eq414}
    	\dot{m} &= (A_1 + \partial_{\bar{n}} I) \dot{n} + \left( s'(r) - s_0 \right) \left( \frac{\lambda_1 \dot{m}}{1 - \frac{\lambda_1}{b_0}} - (A_1 + \partial_{\bar{n}} I) \frac{\lambda_2 \dot{n}}{1 - \frac{\lambda_2}{b_0}} \right) \notag \\
    	&\quad + \left( s_0 - \frac{1}{b_0} \right) \frac{\lambda_1 \dot{m}}{1 - \frac{\lambda_1}{b_0}} - (A_1 + \partial_{\bar{n}} I) \left( s_0 - \frac{1}{b_0} \right) \frac{\lambda_2 \dot{n}}{1 - \frac{\lambda_2}{b_0}}  \\
    	&\quad + (A_1 + \partial_{\bar{n}} I) \frac{s'(r) - \frac{1}{b_0}}{1 - \frac{\lambda_2}{b_0}} \partial_2 \bar{n} - \frac{s'(r) - \frac{1}{b_0}}{1 - \frac{\lambda_1}{b_0}} \partial_1 \bar{m} + \big( A_2 + \partial_\eta I \big) \tau'(r) \notag\quad \text{on } \Gamma_s.
    \end{align}
    By computation, we have
    \begin{equation}\label{eq419}
    	|s'(r)-s_0|\leqslant|\tau(r)|+|r\tau'(r)|.
    \end{equation}

\section{Global solution of the free boundary problem}\label{s5}
\subsection{The construction of local solutions}
According to assumption \textbf{(H1)}, the piston boundary expands uniformly at constant speed $b_0$ inside the spatial interval $r \in [0, r_0]$, implying that perturbation effects emanate from the boundary point
$(r_0, t_0) = \left(r_0, \frac{r_0}{b_0}\right)$. We construct the forward
characteristic cruve $r = \bar{r}(t)$ originating from this point by solving
the Cauchy problem:
\begin{equation*}\label{eq51}
	\begin{cases}
		\frac{d\bar{r}(t)}{dt} = \hat{u}\Big(\frac{t}{\bar{r}(t)}\Big) + \hat{c}\Big(\frac{t}{\bar{r}(t)}\Big), & t < \frac{r_0}{b_0}, \vspace{1ex} \\
		\bar{r}\big(\frac{r_0}{b_0}\big) = r_0,
	\end{cases}
\end{equation*}
and we can solve for the intersection point $(r_s,t_s)$ of $\bar{r}(t)$ and $t=s_0r$.

\begin{figure}[H]
	\centering
	\begin{minipage}[b]{0.5\textwidth}
		\includegraphics[width=\textwidth]{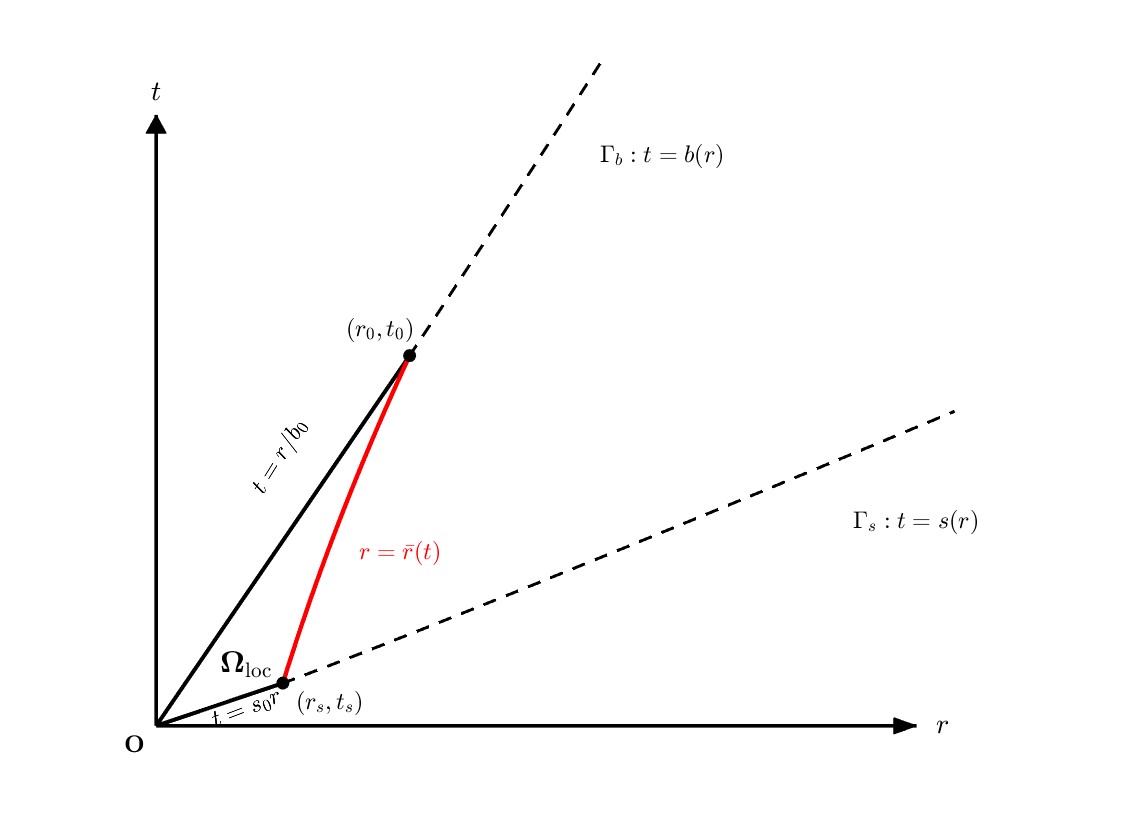}
		\caption{Localized unperturbed domain $\Omega_{\mathrm{loc}}$.}
		\label{Fig0}
	\end{minipage}
\end{figure}

Let $\Omega_{\mathrm{loc}}$ denote the localized space-time domain bounded by the line $r = b_0 t$, $t = s_0 r$, and the characteristic trajectory $r = \bar{r}(t)$, see Figure \ref{Fig0}. 
Inside $\Omega_{\mathrm{loc}}$, the solution to the free boundary problem \eqref{e1}–\eqref{e14} coincides identically with the background self-similar solution, i.e.
\begin{equation*}
	(u, \rho)(r, t) = (\hat{u}, \hat{\rho})\Big(\frac{t}{r}\Big), \quad (r, t) \in \Omega_{\mathrm{loc}}.
\end{equation*}
In order to extend this local solution to global, it is necessary to establish uniform $C^1$ decay estimates for the perturbation variables outside $\Omega_{\mathrm{loc}}$.

\subsection{Uniform $C^1$ a priori estimates}

In this subsection, we establish uniform decay estimates for
both the perturbed variables and the free shock position.

\begin{lemma}\label{le51}
	Let $\Omega_T = \{(r, t) \mid 0< r \leqslant T, \, s(r) \leqslant t \leqslant b(r)\}$ for any large $T$. When $\varepsilon$ is sufficiently small, there exist positive constants $C_0$ and $C_1$ independent of $\varepsilon$ and $T$, the solution of the problem \eqref{e1}, \eqref{eq31} and \eqref{eq32} satisfies
	\begin{equation}\label{e41}
		|\bar{m}| < \frac{C_0 \varepsilon}{r^{\delta_0}}, \quad |\bar{n}| < \frac{C_0 \varepsilon}{r^{\delta_0}}, \quad |\dot{m}| < \frac{C_1 \varepsilon}{r^{1+{\delta_0}}}, \quad |\dot{n}| < \frac{C_1 \varepsilon}{r^{1+{\delta_0}}},\quad\text{in}\ \Omega_T,
	\end{equation}
	and the shock position perturbation satisfies
	\begin{equation}\label{e42}
		|\tau(r)| < \frac{C_0 \varepsilon}{2b_0^2r^{\delta_0}}, \quad |\tau'(r)| < \frac{3C_0 \varepsilon}{b_0^2r^{1+{\delta_0}}} \quad \text{for } r \in\left(0,T\right].
	\end{equation}
\end{lemma}

\begin{proof}
	The proof employs the method of reflected characteristics within a continuous induction framework.
	For any given point $(r_1, t_1) \in \Omega_T$, let $\gamma_j(z; r_1, t_1)$ ($j=1,2$) denote the backward $j$-th characteristic curve starting at $(r_1, t_1)$, which is defined by
	\begin{equation*}
		\begin{cases}
			\frac{d\gamma_j(z; r_1, t_1)}{dz} = \lambda_j \left( u\left(\gamma_j(z; r_1, t_1), z\right) \right), \quad z \leqslant t_1; \\
			\gamma_j(t_1; r_1, t_1) = r_1.
		\end{cases}
	\end{equation*}
	
    As shown in Figure \ref{Fig2}, the backward characteristic curve
    $r=\gamma_1(z; r_1, t_1)$ intersects $\Gamma_s$ at
    $(\phi_1(r_1, t_1), \psi_1(r_1, t_1))$, then tracing the backward $C_2$ characteristic issuing from $(\phi_1, \psi_1)$ to $\Gamma_b$ yields the intersection $(\pi_2(r_1, t_1), \zeta_2(r_1, t_1))$. Symmetrically, the $C_2$ curve $r = \gamma_2(z; r_1, t_1)$ meets $\Gamma_b$ at $(\phi_2(r_1, t_1), \psi_2(r_1, t_1))$, and the corresponding backward $C_1$ path meets $\Gamma_s$ at $(\pi_1(r_1, t_1), \zeta_1(r_1, t_1))$.
		\begin{figure}[H]
		\centering
		\begin{minipage}[b]{0.5\textwidth}
			\includegraphics[width=\textwidth]{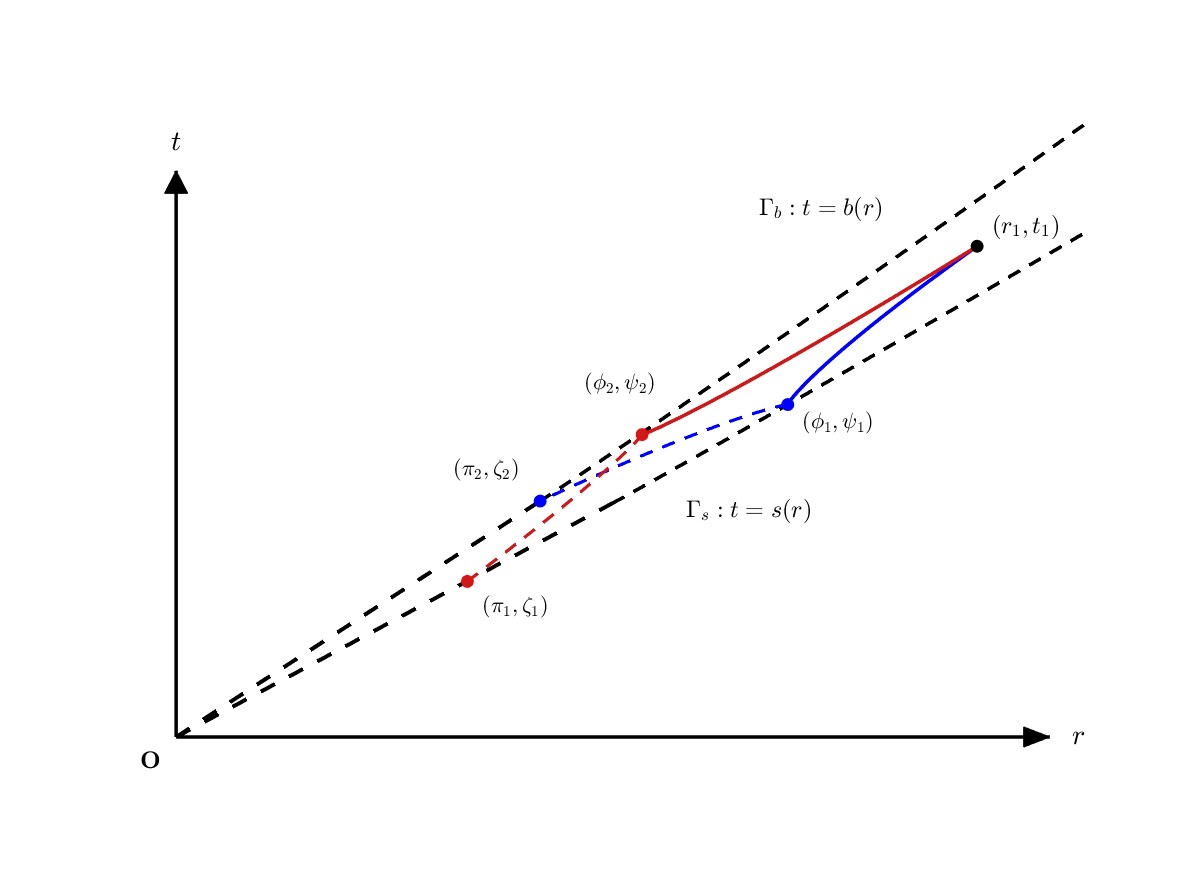}
			\caption{Reflected characteristics curves in $\Omega_T$.}
			\label{Fig2}
		\end{minipage}
	\end{figure}
	
    Let $\Omega(r_1,t_1)$ be a closed domain encircled by $\Gamma_s$, $\Gamma_b$ and $\gamma_j(z; r_1, t_1)$ ($j=1,2$), we need to prove that under the assumption $(\mathcal{H})$
    \[
    (\mathcal{H}): \quad
    \begin{cases}
    	|\bar{m}|,|\bar{n}| < \frac{C_0 \varepsilon}{r^{\delta_0}},\quad \text{in}\ \Omega(r_1,t_1)\backslash\{(r_1,t_1)\}, \\
    	|\dot{m}|,|\dot{n}| < \frac{C_1 \varepsilon}{r^{1+{\delta_0}}},\quad \text{in}\ \Omega(r_1,t_1)\backslash\{(r_1,t_1)\}, \\
    	|\tau(r)| < \frac{C_0 \varepsilon}{2b_0^2r^{\delta_0}}, \quad |\tau'(r)| < \frac{3C_0 \varepsilon}{b_0^2r^{1+{\delta_0}}},\quad \text{on}\ r\in\left(0,\phi_1(r_1,t_1)\right),
    \end{cases}
    \]
     $(\mathcal{\hat{H}})$ is valid
    	\[
    (\hat{\mathcal{H}}): \quad
    \begin{cases}
    	|\bar{m}(r_1,t_1)|,|\bar{n}(r_1,t_1)| < \frac{C_0 \varepsilon}{r^{\delta_0}}, \\
    	|\dot{m}(r_1,t_1)|,|\dot{n}(r_1,t_1)| < \frac{C_1 \varepsilon}{r^{1+{\delta_0}}} \\
    	|\tau(\phi_1(r_1,t_1))| < \frac{C_0 \varepsilon}{2b_0^2r^{\delta_0}}, \quad |\tau'(\phi_1(r_1,t_1))| < \frac{3C_0 \varepsilon}{b_0^2r^{1+{\delta_0}}}.
    \end{cases}
    \]
    
    In the following discussion, let $\epsilon>0$ be an arbitrary small constant.
	We begin by controlling the shock perturbation $\tau(r)$. Via \eqref{eq213}, \eqref{e34} and \eqref{e320}, we get
	\begin{equation}\label{eq55}
		r\tau'+(3+\epsilon)\tau>-(1+\epsilon)\cdot\frac{\bar{m}+\bar{n}}{2b_0^2} \quad \text{and}\quad r\tau'+(3-\epsilon)\tau<(-1+\epsilon)\cdot\frac{\bar{m}+\bar{n}}{2b_0^2}.
	\end{equation}
	A direct computation shows
	\begin{align}\label{eq56}
		|\tau(\phi_1)|&<e^{-\int_{r_s}^{\phi_1}\frac{3-\epsilon}{r}dr}\left(\int_{r_s}^{\phi_1}\left((1+\epsilon)\frac{|\bar{m}|+|\bar{n}|}{2b_0^2r}\cdot e^{\int_{r_s}^{r}\frac{3-\epsilon}{r}dr}\right)dr\right)\notag\\
		&<e^{-\int_{r_s}^{\phi_1}\frac{3-\epsilon}{r}dr}\left(\int_{r_s}^{\phi_1}\left((1+\epsilon)\frac{C_0\varepsilon}{b_0^2r^{1+\delta_0}}\cdot e^{\int_{r_s}^{r}\frac{3-\epsilon}{r}dr}\right)dr\right)\\
		&=\left(\frac{1+\epsilon}{3-\epsilon-\delta_0}\right)\cdot\left[\frac{C_0\varepsilon}{b_0^2\phi_1^{\delta_0}}-\frac{C_0\varepsilon}{b_0^2\phi_1^{3-\epsilon}}r_s^{3-\epsilon-\delta_0}\right]\notag\\
		&<\frac{C_0\varepsilon}{2b_0^2\phi_1^{\delta_0}},\notag
	\end{align}
	subsituting \eqref{eq56} into (\ref{eq55}), we have
	\begin{equation*}\label{eq57}
		|\tau'(\phi_1)|<\frac{3C_0\varepsilon}{b_0^2\phi_1^{1+\delta_0}}.
	\end{equation*}

	Next, we establish the decay bounds for $|\bar{m}|$ and $|\bar{n}|$. When $\varepsilon$ and $\rho_0$ are both sufficiently small, there holds
	\begin{align}
		&\left\lvert\frac{F(m,n)-F(\hat{m},\hat{n})}{r\lambda_1(m,n)}\right\rvert+\left\lvert\frac{(\lambda_1(\hat{m},\hat{n})-\lambda_1(m,n))\partial_r\hat{m}}{\lambda_1(m,n)}\right\rvert<\frac{\gamma+\epsilon}{r}(|\bar{m}|+|\bar{n}|),\label{eq581}\\
		&\left\lvert-\frac{F(m,n)-F(\hat{m},\hat{n})}{r\lambda_2(m,n)}\right\rvert+\left\lvert\frac{(\lambda_2(\hat{m},\hat{n})-\lambda_2(m,n))\partial_r\hat{n}}{\lambda_2(m,n)}\right\rvert<\frac{\gamma+\epsilon}{r}(|\bar{m}|+|\bar{n}|),\label{eq582}
	\end{align}
	with \eqref{e313}.
	Besides, \eqref{e214} and \eqref{eq 213} show
	\begin{equation}\label{eq59}
	\frac{\phi_j(r_1,t_1)}{r_1}>1-\frac{1}{(3-2\epsilon)\sqrt{\gamma}}b_0^{-\frac{1}{\gamma}}\rho_0^\frac{\gamma-1}{2\gamma},\quad	 \frac{\pi_j(r_1,t_1)}{r_1}>1-\frac{2}{(3-2\epsilon)\sqrt{\gamma}}b_0^{-\frac{1}{\gamma}}\rho_0^\frac{\gamma-1}{2\gamma},\ j=1,2.
	\end{equation}
	
	Integrating \eqref{eq37} along characteristic curves, we get 
	\begin{align*}
		|\bar{m}(r_1,t_1)|&\leqslant|\bar{m}(\phi_1,\psi_1)|+(2\gamma+\epsilon)\int_{\phi_1}^{r_1}\frac{C_0\varepsilon}{r^{1+\delta_0}}dr\notag\\
		&\leqslant|A_1\bar{n}(\phi_1,\psi_1)|+|A_2\tau(\phi_1)|+|I(\bar{n}(\phi_1,\psi_1),\tau(\phi_1))|+(2\gamma+\epsilon)\int_{\phi_1}^{r_1}\frac{C_0\varepsilon}{r^{1+\delta_0}}dr\\
		&\leqslant|A_1\bar{n}(\pi_2,\zeta_2)|+|A_2\tau(\phi_1)|+|I(\bar{n}(\phi_1,\psi_1),\tau(\phi_1))|+(2\gamma+\epsilon)\int_{\pi_2}^{r_1}\frac{C_0\varepsilon}{r^{1+\delta_0}}dr\notag\\
		&\leqslant|A_1\bar{m}(\pi_2,\zeta_2)|+|G_{11}(\pi_2)|+|G_{12}(\pi_2,\zeta_2)|+|A_2\tau(\phi_1)|+|I(\bar{n}(\phi_1,\psi_1),\tau(\phi_1))|\\
		&+(2\gamma+\epsilon)\int_{\pi_2}^{r_1}\frac{C_0\varepsilon}{r^{1+\delta_0}}dr,\notag
	\end{align*}
	with \eqref{eq316}, \eqref{eq320}, \eqref{eq581} and \eqref{eq582}, which means
	\begin{align}\label{e48}
		|\bar{m}(r_1,t_1)|&\leqslant\frac{1 -  \frac{(4-\epsilon)b_0^{-\frac{1}{\gamma}}}{\sqrt{\gamma}} \rho_0^\frac{\gamma-1}{2\gamma}}{\left(1-\frac{2}{(3-2\epsilon)\sqrt{\gamma}}b_0^{-\frac{1}{\gamma}}\rho_0^\frac{\gamma-1}{2\gamma}\right)^{\delta_0}}\cdot\frac{C_0\varepsilon}{r_1^{\delta_0}}+\frac{4b_0^2+\epsilon}{\left(1-\frac{2}{(3-2\epsilon)\sqrt{\gamma}}b_0^{-\frac{1}{\gamma}}\rho_0^\frac{\gamma-1}{2\gamma}\right)^{\delta_0}}\cdot\frac{\varepsilon}{r_1^{\delta_0}}+\frac{\frac{5b_0^2}{1-\delta_0}+\epsilon}{\left(1-\frac{2}{(3-2\epsilon)\sqrt{\gamma}}b_0^{-\frac{1}{\gamma}}\rho_0^\frac{\gamma-1}{2\gamma}\right)^{\delta_0}}\cdot\frac{\varepsilon}{r_1^{\delta_0}}\notag\\
		&+\frac{\frac{(2+\epsilon) b_0^{-\frac{1}{\gamma}}}{\sqrt{\gamma}}\rho_0^\frac{\gamma-1}{2\gamma}}{\left(1-\frac{1}{(3-2\epsilon)\sqrt{\gamma}}b_0^{-\frac{1}{\gamma}}\rho_0^\frac{\gamma-1}{2\gamma}\right)^{\delta_0}}\cdot\frac{C_0\varepsilon}{r_1^{\delta_0}}+\mathcal{M}\left[\left(\frac{C_0\varepsilon}{r_1^{\delta_0}}\right)^2+\left(\frac{C_0\varepsilon}{2b_0^2r_1^{\delta_0}}\right)^2\right]\notag\\
		&+\frac{2\gamma+\epsilon}{\delta_0}\cdot\left[\frac{1}{\left(1-\frac{2}{(3-2\epsilon)\sqrt{\gamma}}b_0^{-\frac{1}{\gamma}}\rho_0^\frac{\gamma-1}{2\gamma}\right)^{\delta_0}}-1\right]\cdot\frac{C_0\varepsilon}{r_1^{\delta_0}}\\
		&<\left[1-
		\left(1-\frac{2\gamma+\delta_0+\epsilon}{3-2\epsilon}-\epsilon\right)
		\cdot\frac{2}{\sqrt{\gamma}}b_0^{-\frac{1}{\gamma}}\rho_0^\frac{\gamma-1}{2\gamma}+\mathcal{D}\rho_0^\frac{\gamma-1}{\gamma}+\mathcal{M}\left(\frac{C_0\varepsilon}{r_1^{\delta_0}}+\frac{C_0\varepsilon}{4b_0^4r_1^{\delta_0}}\right)\right]\cdot\frac{C_0\varepsilon}{r_1^{\delta_0}}\notag\\
		&+\left[\left(1+\frac{2\delta_0}{(3-2\epsilon)\sqrt{\gamma}}b_0^{-\frac{1}{\gamma}}\rho_0^\frac{\gamma-1}{2\gamma}\right)\cdot\left(4b_0^2+\frac{5b_0^2}{1-\delta_0}+2\epsilon\right)\right]\cdot\frac{\varepsilon}{r_1^{\delta_0}},\notag
	\end{align}
	with \eqref{eq3161}, \eqref{3162}, \eqref{eq:3.22} and \eqref{eq581}-\eqref{eq59} under assuption $(\mathcal{H})$, where $\mathcal{D}$ is a constant depending on $b_0$, $\delta_0$ and $\gamma$. It is obvious that the bound of $|\bar{m}(r_1,t_1)|$ in \eqref{e48} is a continuous function associated to $\epsilon$, it holds that
	\begin{equation}
			\lim\limits_{\epsilon\to 0^+}1-\frac{2\gamma+\delta_0+\epsilon}{3-2\epsilon}-\epsilon>0,
	\end{equation}
	due to $\delta_0<3-2\gamma$. By choosing $C_0=\frac{9\sqrt{\gamma}b_0^{2+\frac{1}{\gamma}}}{(1-\frac{2\gamma+\delta_0}{3})\rho_0^\frac{\gamma-1}{2\gamma}}$, we can get $|\bar{m}(r_1, t_1)| <\frac{C_0\varepsilon}{r_1^{\delta_0}}$.
	
	Similarly, $|\bar{n}(r_1, t_1)| <\frac{C_0\varepsilon}{r_1^{\delta_0}}$ holds.

	Finally, we close the estimates for $|\dot{m}|$ and
	$|\dot{n}|$. From the definition of $\partial_j$ ($j=1, 2$), \eqref{e233} and \eqref{eq581}-\eqref{eq582}, we have
	\begin{equation}\label{eq410}
		|\partial_1\bar{m}|,|\partial_2\bar{n}|<\frac{\gamma b_0+\epsilon}{r}(|\bar{m}|+|\bar{n}|).
	\end{equation}
	Thus, using \eqref{eq212}-\eqref{eq215}, \eqref{e233}, \eqref{eq43} , \eqref{eq45}-\eqref{e337} and \eqref{eq410}, we have
	 \begin{align*}
		\left|\frac{1}{r} \partial_p \big[ F(m, n) - F(\hat{m}, \hat{n}) \big]\right|&<\left[\frac{(\gamma-1)b_0}{2}+\epsilon\right]\cdot\frac{|\dot{m}|+|\dot{n}|}{r},\\
		\left|\partial_p \big[ \lambda_1(\hat{m}, \hat{n}) - \lambda_1(m, n) \big](\partial_r \hat{m} + \partial_r \bar{m})\right|&<\frac{b_0+\epsilon}{2}\cdot\frac{(1+\gamma)|\dot{m}|+(3-\gamma)|\dot{n}|}{r},\\
		\left|\big[ \lambda_1(\hat{m}, \hat{n}) - \lambda_1(m, n) \big]\partial_p \partial_r \hat{m}\right|&<\frac{b_0+\epsilon}{2}\cdot\frac{(1+\gamma)|\bar{m}|+(3-\gamma)|\bar{n}|}{r^2},\\
		\left|\partial_p \lambda_1(\hat{m}, \hat{n}) \partial_r \bar{m}\right|&<\epsilon\left(\frac{|\bar{m}|}{r^2}+\frac{|\dot{m}|}{r}\right),\\
		\left|\frac{F(m, n) - F(\hat{m}, \hat{n})}{r^2}\right|&<\left[\frac{(\gamma-1)b_0}{2}+\epsilon\right]\cdot\frac{|\bar{m}|+|\bar{n}|}{r^2}.
	\end{align*}
	The same argument can apply to $H_2$. Thus, we have

	\begin{equation}\label{eq511}
		\left\lvert\frac{H_i(m,n)}{\lambda_i(m,n)}\right\rvert<(1+\epsilon)\left[\frac{\gamma|\bar{m}|+|\bar{n}|}{r^2}+\frac{\gamma|\dot{m}|+|\dot{n}|}{r}\right],\quad \text{in}\ \Omega(r_1,t_1),
	\end{equation}
	where $i=1,2$.
	
	For the fixed boundary condition \eqref{eq:4.12} along $\Gamma_b$, combining with (\textbf{H1})-(\textbf{H2}), \eqref{eq212}-\eqref{eq215}, \eqref{e233}, \eqref{eq416} and \eqref{eq410}, we have
	\begin{equation*}
		|\dot{n}|<|\dot{m}|+\frac{\mathcal{C}'}{\rho_0^\frac{\gamma-1}{2\gamma}}\cdot\left(|\dot{m}|+|\dot{n}|\right)\cdot\frac{\varepsilon}{r^{\delta_0}}+\frac{\mathcal{C}'}{\rho_0^\frac{\gamma-1}{2\gamma}}\cdot\left(|\bar{m}|+|\bar{n}|\right)\cdot\frac{\varepsilon}{r^{1+\delta_0}}+\mathcal{C}'\frac{\varepsilon}{r^{1+\delta_0}}+\left|\frac{\frac{2\varepsilon}{r^{1+\delta_0}}}{1-\left(\frac{1}{b_0}+\frac{\varepsilon}{r^{\delta_0}}\right)^2}\right|,
	\end{equation*}
	where $\mathcal{C}'$ depends on $\delta_0$ and $b_0$. Denote $\mathcal{L}_b$ is a constant that depends on $\delta_0$ and $b_0$, when $\varepsilon$ is sufficiently small, it yields
	\begin{equation}\label{eq 48}
		|\dot{n}|<(1+\epsilon)|\dot{m}|+\mathcal{L}_b\cdot\frac{\varepsilon}{r^{1+\delta_0}}.
	\end{equation}

	For the free boundary condition (\ref{eq414}) along $\Gamma_s$, when $\varepsilon$ is sufficiently small, $A_1 + \partial_{\bar{n}}I$ amd $A_2 + \partial_\eta I$ satisfy
	\begin{align*}
		&|A_1 + \partial_{\bar{n}}I|< 1 -\frac{(4-\epsilon)b_0^{-\frac{1}{\gamma}}}{\sqrt{\gamma}}\rho_0^\frac{\gamma-1}{2\gamma} +\mathcal{M}(|\bar{n}|+|\tau|)<1 -\frac{(4-2\epsilon)b_0^{-\frac{1}{\gamma}}}{\sqrt{\gamma}}\rho_0^\frac{\gamma-1}{2\gamma},\\
		  &|A_2 + \partial_\eta I|<\frac{4+\epsilon}{\sqrt{\gamma}}b_0^{2-\frac{1}{\gamma}} \rho_0^\frac{\gamma-1}{2\gamma}+\mathcal{M}(|\bar{n}|+|\tau|)<\frac{4+2\epsilon}{\sqrt{\gamma}}b_0^{2-\frac{1}{\gamma}} \rho_0^\frac{\gamma-1}{2\gamma},
	\end{align*}
	with \eqref{eq314}. And there holds
	\begin{align*}
		|\dot{m}|<&|A_1 + \partial_{\bar{n}}I|\cdot|\dot{n}|+\frac{C_0+\mathcal{C}'}{\rho_0^\frac{\gamma-1}{2\gamma}}\cdot\left(|\dot{m}|+|A_1 + \partial_{\bar{n}}I|\cdot|\dot{n}|\right)\cdot\frac{\varepsilon}{r^\delta_0}\\
		&+\frac{C_0+\mathcal{C}'}{\rho_0^\frac{\gamma-1}{2\gamma}}\cdot\left(1+|A_1 + \partial_{\bar{n}}I|\right)\cdot\left(|\bar{m}|+|\bar{n}|\right)\cdot\frac{\varepsilon}{r^{1+\delta_0}}\\
		&+|A_2 + \partial_\eta I|\cdot|\tau'(r)|,
	\end{align*}
	from \eqref{eq414} with \eqref{e216}, \eqref{eq212}, \eqref{eq214}, \eqref{e233}, \eqref{eq419} and \eqref{eq410}. Thus
	\begin{equation}\label{eq521}
		|\dot{m}|<(1+\epsilon)\cdot|A_1 +\partial_{\bar{n}}I|\cdot|\dot{n}|+|A_2 + \partial_\eta I|\cdot|\tau'|,\quad \text{as}\ (r,t)\in\left\{t=s(r), r\leqslant\phi_1(r_1,t_1)\right\}.
	\end{equation}
	
	Integrating \eqref{eq44}, together with \eqref{eq511}-\eqref{eq521} and assuption $(\mathcal{H})$, we get
	\begin{align*}
		|\dot{m}(r_1,t_1)|\leqslant&|\dot{m}(\phi_1,\psi_1)|+(\gamma+1+\epsilon)\int_{\phi_1}^{r_1}\left(\frac{C_0\varepsilon}{r^{2+\delta_0}}+\frac{C_1\varepsilon}{r^{2+\delta_0}}\right)dr\\
		\leqslant&(1+\epsilon)\cdot|A_1 + \partial_{\bar{n}}I|\cdot|\dot{n}(\phi_1,\psi_1)|+|A_2 + \partial_\eta I|\cdot|\tau'(\phi_1)|+(\gamma+1+\epsilon)\int_{\phi_1}^{r_1}\left(\frac{C_0\varepsilon}{r^{2+\delta_0}}+\frac{C_1\varepsilon}{r^{2+\delta_0}}\right)dr\notag\\
		\leqslant&(1+\epsilon)\cdot|A_1 + \partial_{\bar{n}}I|\cdot|\dot{n}(\pi_2,\zeta_2)|+|A_2 + \partial_\eta I|\cdot|\tau'(\phi_1)|+(\gamma+1+\epsilon)\int_{\pi_2}^{r_1}\left(\frac{C_0\varepsilon}{r^{2+\delta_0}}+\frac{C_1\varepsilon}{r^{2+\delta_0}}\right)dr\notag\\
		\leqslant&(1+\epsilon)^2\cdot|A_1 + \partial_{\bar{n}}I|\cdot|\dot{m}(\pi_2,\zeta_2)|+\frac{\mathcal{L}_b\varepsilon}{\pi_2^{1+\delta_0}}+|A_2 + \partial_\eta I|\cdot|\tau'(\phi_1)|\\
		&+(\gamma+1+\epsilon)\int_{\pi_2}^{r_1}\left(\frac{C_0\varepsilon}{r^{2+\delta_0}}+\frac{C_1\varepsilon}{r^{2+\delta_0}}\right)dr,\notag
	\end{align*}
	based on (\ref{eq59}), we get
	\begin{align*}
		|\dot{m}(r_1,t_1)|\leqslant&(1+\epsilon)^2\cdot\frac{1 -  \frac{(4-2\epsilon)b_0^{-\frac{1}{\gamma}}}{\sqrt{\gamma}} \rho_0^\frac{\gamma-1}{2\gamma}}{\left(1-\frac{2}{(3-2\epsilon)\sqrt{\gamma}}b_0^{-\frac{1}{\gamma}}\rho_0^\frac{\gamma-1}{2\gamma}\right)^{1+\delta_0}}\cdot\frac{C_1\varepsilon}{r_1^{1+\delta_0}}+\frac{\mathcal{L}_b}{\left(1-\frac{2}{(3-2\epsilon)\sqrt{\gamma}}b_0^{-\frac{1}{\gamma}}\rho_0^\frac{\gamma-1}{2\gamma}\right)^{1+\delta_0}}\cdot\frac{\varepsilon}{r_1^{1+\delta_0}}\\\notag
		&+\frac{3\cdot\frac{4+2\epsilon}{\sqrt{\gamma}}b_0^{-\frac{1}{\gamma}} \rho_0^\frac{\gamma-1}{2\gamma}}{\left(1-\frac{2}{(3-2\epsilon)\sqrt{\gamma}}b_0^{-\frac{1}{\gamma}}\rho_0^\frac{\gamma-1}{2\gamma}\right)^{1+\delta_0}}\cdot\frac{C_0\varepsilon}{r_1^{1+\delta_0}}
		+\frac{\gamma+1+\epsilon}{1+\delta_0}\left[\frac{1}{\left(1-\frac{2}{(3-2\epsilon)\sqrt{\gamma}}b_0^{-\frac{1}{\gamma}}\rho_0^\frac{\gamma-1}{2\gamma}\right)^{1+\delta_0}}-1\right]\frac{(C_0+C_1)\varepsilon}{r_1^{1+\delta_0}}\\\notag
		<&\left\{(1+\epsilon)^2-\left[(1+\epsilon)^2\cdot  \frac{(4-2\epsilon)b_0^{-\frac{1}{\gamma}}}{\sqrt{\gamma}}-(1+\epsilon)^2\cdot\frac{2(1+\delta_0)b_0^{-\frac{1}{\gamma}}}{(3-2\epsilon)\sqrt{\gamma}}-\frac{2(\gamma+1+\epsilon)b_0^{-\frac{1}{\gamma}}}{(3-2\epsilon)\sqrt{\gamma}}\right]\rho_0^\frac{\gamma-1}{2\gamma}+\mathcal{D}\rho_0^\frac{\gamma-1}{\gamma}\right\}\cdot\frac{C_1\varepsilon}{r_1^{1+\delta_0}}\\
		&+\left\{\mathcal{L}_b+\left[\frac{2\mathcal{L}_b(1+\delta_0)}{3-2\epsilon}+3(4+2\epsilon)C_0+\frac{2C_0(\gamma+1+\epsilon)}{3-2\epsilon}\right]\frac{b_0^{-\frac{1}{\gamma}}}{\sqrt{\gamma}}\rho_0^\frac{\gamma-1}{2\gamma}+\mathcal{D}\rho_0^\frac{\gamma-1}{\gamma}\right\}\frac{\varepsilon}{r_1^{1+\delta_0}},
	\end{align*}
	by choosing suitable $C_1$, we have $|\dot{m}(r_1,t_1)|<\frac{C_1\varepsilon}{r_1^{1+\delta_0}}$. Similarly, $|\dot{n}(r_1, t_1)| <\frac{C_1\varepsilon}{r_1^{1+\delta_0}}$ holds.
\end{proof}

\subsection{Proof of Theorem \ref{th1}}
 By \eqref{e233}, \eqref{n36}, \eqref{n38}, \eqref{eq36}, \eqref{eq43},  \eqref{eq419} and Lemma \ref{le51}, we can obtain that estimates \eqref{e18}-\eqref{e11} hold. Thus, the global existence and unqueness of shock solutions are established by the local existence result and step-by-step continuation arguments for quasilinear hyperbolic free boundary problems.

\section{Acknowledgment}
Huimin Yu is supported by NSFC Grant No. 12271310 and Natural Science Foundation of Shandong
Province ZR2022MA088.

\end{document}